%% file: main.tex
\documentclass[a4paper,10pt]{amsart}
\usepackage[USenglish]{babel}
\usepackage[utf8]{inputenc}
\pdfoutput=1 

\usepackage[style=american]{csquotes}

\usepackage{amsmath,amsthm,amscd,amssymb}
\usepackage{helvet}
\usepackage{courier}
\usepackage{graphicx}
\usepackage{multicol}
\usepackage{footmisc}
\usepackage{tikz}
\usetikzlibrary{calc}
\usepackage{color}
\usepackage{hyperref}
\usepackage{url}
\usepackage{xspace}
\usepackage{listings}

\definecolor{appcolor}{rgb}{0,0,0}
\definecolor{commentcolor}{rgb}{.3,.3,.3}
\definecolor{stringcolor}{RGB}{103,5,172}
\definecolor{typecolor}{RGB}{140,75,0}
\definecolor{propcolor}{RGB}{0,100,14}
\definecolor{light_gray}{gray}{0.7}

\newcommand{\argmax}{\textnormal{arg max}}
\newcommand{\argmin}{\textnormal{arg min}}
\newcommand{\curly}[1]{\mathcal{#1}}
\DeclareMathOperator{\conv}{conv}
\DeclareMathOperator{\seccone}{seccone}

\newcommand{\atint}{\texttt{a-tint}\xspace}
\newcommand{\polymake}{\texttt{poly\-make}\xspace}
\newcommand{\Singular}{\texttt{Singular}\xspace}
\newcommand{\Macaulay}{\texttt{Macaulay2}\xspace}
\newcommand{\shelltrop}{{\color{appcolor}tropical $>$}}
\newcommand{\shellfan}{{\color{appcolor}fan $>$}}
\newcommand{\shellideal}{{\color{appcolor}ideal $>$}}
\newcommand{\shellpoly}{{\color{appcolor}polytope $>$}}

\newcommand\1{{\mathbf 1}}
\newcommand{\N}{\mathbb{N}}
\newcommand{\Z}{\mathbb{Z}}
\newcommand{\T}{\mathbb{T}}
\newcommand{\Q}{\mathbb{Q}}
\newcommand{\R}{\mathbb{R}}

\newcommand{\TT}[1]{\R^{#1+1}/\R\1}

\newcommand{\Dr}{\textnormal{Dr}}
\newcommand{\Gr}{\textnormal{Gr}}
\newcommand{\TGr}{\textnormal{TGr}}
\newcommand{\Sym}{\textnormal{Sym}}
\DeclareMathOperator\tdet{tdet}

\theoremstyle{definition}
\newtheorem{definition}{Definition}[section]

\newtheorem*{acknowledgement}{Acknowledgement}

\newtheorem{remark}[definition]{Remark}
\theoremstyle{plain}
\newtheorem{theorem}[definition]{Theorem}

\lstdefinelanguage{pmshell}
{
  basicstyle=\small\ttfamily,
  keywords=[1]{Hypersurface,ValuatedMatroid,SubdivisionOfPoints,Ideal,TropicalNumber,Min,Max,Matrix, Rational,Vector,polytope,matroid,Int,fan,IncidenceMatrix,Array,Set},
  morekeywords=[2]{N_ISOLATED, N_FAMILIES, POLYNOMIAL, DEGREE, GENUS, BASES, VALUATION_ON_BASES, WEIGHTS, N_ELEMENTS, FAR_VERTICES, CONES, POINTS, VERTICES, TIGHT_SPAN, MAXIMAL_POLYTOPES, MONOMIALS, COEFFICIENTS, MAXIMAL_CELLS, RAYS, N_RAYS, N_MAXIMAL_CELLS, FACETS, LINEAR_SPAN, CURVE_EDGE_LENGTHS, VISUAL, GENERATORS, GROEBNER, ORDER_VECTOR,INITIAL_FORMS,RANK,N_POINTS,N_MAXIMAL_POLYTOPES,TUTTE_POLYNOMIAL, NESTED_COEFFICIENTS},
  morekeywords=[3]{map, grep, join, for, my, secondary_cone,tdet,tdet_and_perm,elem,dual_subdivision,toTropicalPolynomial,minor,transpose,row,size,ones_vector,labeled,hypersimplex,uniform_matroid,rows,check_basis_exchange_axiom,pluecker_ideal,contains_monomial,linear_space,fano_matroid,matroid_fan,trivial_valuation,matroid_fan_from_flats,intersect,intersect_in_smooth_surface,lines_in_cubic, representative,all_families,visualize_in_surface,matroid_ring_cycle,LengthLabels,nested_matroids,direct_sum},
  stringstyle=\color{stringcolor},
  keywordstyle=[1]\color{typecolor},
  keywordstyle=[2]\color{propcolor},
  keywordstyle=[3]\color{typecolor},
  numbers=none,
  captionpos=b,
  showspaces=false,
  showstringspaces=false,
  morestring=[b]",
  escapechar=&,
  frame=single,
}

\begin{document}

\title{Tropical Computations in polymake}
\author{Simon Hampe and Michael Joswig}
\address{Technische Universität Berlin\\
Institut für Mathematik, Sekretariat MA 6-2\\
Straße des 17. Juni 136, 10623 Berlin}
\email{(hampe,joswig)@math.tu-berlin.de} 

\begin{abstract}
We give an overview of recently implemented \polymake features for computations in tropical geometry.
The main focus is on explicit examples rather than technical explanations.
Our computations employ tropical hypersurfaces, moduli of tropical plane curves, tropical linear spaces and Grassmannians, lines on tropical cubic surfaces as well as intersection rings of matroids.
\end{abstract}

\maketitle

\section{Introduction}

Many avenues lead to tropical geometry as we know it today.
One motivation comes from studying algebraic varieties (over some field with a non-Archimedean valuation) via their piecewise-linear images (under the valuation map).
This is useful since many interesting properties are preserved, and they often become algorithmically accessible via tools from polyhedral geometry \cite{Tropical+Book}.

Many of these methods are actually implemented, and we start out with giving a brief overview.
The standard software to compute tropical varieties with constant coefficients is Jensen's \texttt{Gfan} \cite{gfan}.
Its main function is to traverse the dual graph of the Gr\"obner fan of an ideal and to construct the associated tropical variety as a subfan.
The \texttt{Singular} \cite{singular} library \texttt{tropical.lib} \cite{tropical.lib} by Jensen, Markwig, Markwig and Ren interfaces to \texttt{Gfan} and implements extra functionality on top.
The most recent version also covers Ren's implementation of tropical varieties with arbitrary coefficients \cite{Ren:2015}.
Rinc\'on's program \texttt{TropLi} computes tropical linear spaces from matrix input \cite{rinconcompute}, whereas the \texttt{Tropical Polyhedra Library} by Allamigeon allows to manipulate tropical polyhedra \cite{tplib}.
The \polymake system is a comprehensive system for polyhedral geometry and adjacent areas of discrete mathematics.
Basic support for computations with tropical hypersurfaces and tropical polytopes goes back as far as version 2.0 from 2004.
A much more substantial contribution to \polymake was the extension \texttt{a-tint} for tropical intersection theory \cite{hatint}.

In this paper we report on the recent rewrite of all functionality related to tropical geometry in \polymake.
This largely builds on \texttt{a-tint}, which is now a bundled extension and which itself has undergone a massive refactoring.
Here we refer to the current version 3.0 of \polymake from 2016.

Our paper is organized as follows. We start out with the basics of tropical arithmetic and tropical matrix operations. This topic connects tropical geometry to combinatorial optimization \cite{Schrijver03:CO_A}. The most basic geometric objects in our investigation are tropical hypersurfaces. These are the vanishing loci of tropical polynomials. Since the latter are equivalent to finite point sets in $\Z^d$, equipped with real-valued lifting functions, their study is closely related to regular subdivisions \cite{Triangulations}. An interesting new vein in tropical geometry are applications to economics. As an example, we look at arrangements of tropical hypersurfaces as they occur in the product-mix auctions of Baldwin and Klemperer \cite{Klemperer:2009, analysedemand}.
For a regular subdivision $\Sigma$ of a point configuration $P$, the secondary cone comprises all lifting functions, which induce $\Sigma$. 
That cone forms a stratum in the moduli space of tropical hypersurfaces with support set $P$. We exhibit an example computation concerning moduli of tropical plane curves of genus three \cite{BJMS:2015}.

Going from hypersurfaces to more general tropical varieties is a major step. Historically, the first explicit computations dealt with the tropicalization of the Grassmannians \cite{SpeyerSturmfels04}. The classical Grassmannians are the moduli spaces of linear subspaces in a complex vector space. Their tropical analogues parametrize those tropical linear spaces, which arise as tropicalizations. We explore the combinatorics of one tropical linear space. Employing the \polymake interface to \Singular, we verify that it is realizable. Tropical linear spaces are interesting in their own right for their connection with matroid theory \cite{Kapranov93,speyer}. We briefly compare several polyhedral structures on the Bergman fan of a matroid.

A famous classical result of Cayley and Salmon states that every smooth cubic  surface in $\mathbb{P}^3$ over an algebraically closed field contains exactly 27 lines. Vigeland studied the question, whether a similar result holds in the tropical setting \cite{Vigeland}. Based on a \polymake computation with \atint's specialized algorithms for tropical intersection theory, we exhibit a generic tropical cubic surface, $V$, which does not match any of the types listed by Vigeland. The surface $V$ contains 26 isolated lines and three infinite families.
Cohomological methods are an indispensable tool in modern algebraic geometry.
Tropical intersection theory is a first step towards a similar approach to tropical varieties \cite{mikhalkinenumerative, AllermannRau:2010}.
Interestingly, tropical intersection theory is also useful in combinatorics, if applied to tropical linear spaces associated with matroids.
Our final example computation shows how the Tutte polynomial of a matroid can be computed from the nested components \cite{hampematroidring}.

In addition to the features presented here, \polymake also provides functions for computing with Puiseux fractions \cite{JoswigLohoLorenzSchroeter:2016}, tropical polytopes \cite{DevelinSturmfels04,Tropical+Convex+Hull+Computations}, general tropical cycles, tropical morphisms and rational functions \cite{hatint}.

\section{Arithmetic and Linear Algebra}

The \emph{tropical semiring} $\T$ is the set $\R\cup\{\infty\}$ equipped with $\oplus:=\min$ as the \emph{tropical addition} and $\odot:=+$ as the \emph{tropical multiplication}.
Clearly one could also use $\max$ instead of $\min$, but here we will stick to $\oplus=\min$.
In \polymake there is a corresponding data type, which allows to compute in $\Q\cup\{\infty\}$, e.g., the following.

\begin{lstlisting}[language=pmshell,caption={Adding and multiplying tropically.},label=lst:arith]
&\shellpoly& application "tropical";
&\shelltrop& $a = new TropicalNumber<Min>(3);
&\shelltrop& $b = new TropicalNumber<Min>(5);
&\shelltrop& $c = new TropicalNumber<Min>(8);
&\shelltrop& print $a*$c, ", ", $b*$c;
11, 13
&\shelltrop& print (($a + $b) * $c);
11
&\shelltrop& print $a * (new TropicalNumber<Min>("inf"));
inf
\end{lstlisting}
Note that polymake is organized into several \emph{applications}, which serve to separate the various functionalities. Most of our computations take place in the application \texttt{tropical}, but we will occasionally make use of other types of objects, such as matroids, fans and ideals. One can either switch to an application as shown in Listing \ref{lst:arith} or prefix the corresponding types and commands with the name of the application and two colons, such as \texttt{matroid::}. We will see examples below. The software system \polymake is a hybrid design, written in C++ and Perl.
In the \polymake shell the user's commands are interpreted in an enriched dialect of Perl.
Note that here the usual operators ``\texttt{+}'' and ``\texttt{*}'' are overloaded, i.e., they are interpreted as tropical matrix addition and tropical matrix multiplication, respectively.
It is always necessary to explicitly specify the tropical addition via the template parameter \texttt{Min} or \texttt{Max}.
Mixing expressions with \texttt{Min} and \texttt{Max} is not defined and results in an error.
Templates are not part of standard Perl but rather part of \polymake's Perl enrichment.

The type \texttt{TropicalNumber} may be used for coefficients of vectors, matrices and polynomials.
Matrix addition and multiplication are defined --- and interpreted tropically.
Here is a basic application of tropical matrix computations:
Let $A=(a_{ij})\in\R^{d\times d}$ be a square matrix encoding edge lengths on the complete directed graph $K_d$.
If there are no directed cycles of negative length then there is a well defined shortest path between any two nodes, which may be of infinite length.
These shortest path lengths are given by the so-called \emph{Kleene star}
\begin{equation}\label{eq:kleene}
  A^* \ := \ I \oplus A \oplus (A\odot A) \oplus (A \odot A \odot A ) \oplus \dots \enspace ,
\end{equation}
where $I$ is the tropical identity matrix, which has zeros on the diagonal and infinity as a coefficient otherwise.
The assumption that there are no directed cycles of negative length makes the above tropical sum of tropical matrix powers stabilize after finitely many steps.
The direct evaluation of (\ref{eq:kleene}) is precisely the Floyd--Warshall-Algorithm, for computing all shortest paths, known from combinatorial optimization \cite[\S8.4]{Schrijver03:CO_A}.
In Listing~\ref{lst:matrix} below we compute the Kleene star of a $3{\times}3$-matrix, called $A$.
Here we verify that $I\oplus A = I \oplus A \oplus (A\odot A)$, which implies $A^*=I\oplus A$, i.e., all shortest paths are direct.
\begin{lstlisting}[language=pmshell,caption={Adding and multiplying matrices tropically to obtain the Kleene star $A^*$.},label=lst:matrix]
&\shelltrop& $A = new Matrix<TropicalNumber<Min>>(
  [[1,2,3],[1,2,4],[1,0,1]]);
&\shelltrop& $I = new Matrix<TropicalNumber<Min>>(
  [[0,"inf","inf"],["inf",0,"inf"],["inf","inf",0]]);
&\shelltrop& print $I + $A;
0 2 3
1 0 4
1 0 0
&\shelltrop& print $I + $A + $A*$A;
0 2 3
1 0 4
1 0 0
\end{lstlisting}

The \emph{tropical determinant} of $A$ is defined as
\[
\begin{split}
  \tdet A \ :=& \ \bigoplus_{\sigma\in\Sym(d)} a_{1,\sigma(1)} \odot \dots \odot a_{d,\sigma(d)} \\
  =& \ \min\bigl\{a_{1,\sigma(1)}+\dots+a_{d,\sigma(d)}\bigm| \sigma\in\Sym(d)\bigr\} \enspace ,
\end{split}
\]
where $\Sym(d)$ denotes the symmetric group of degree $d$.
This arises from tropicalizing Leibniz' formula for the classical determinant.
Notice that evaluating the tropical determinant is tantamount to solving a linear assignment problem from combinatorial optimization.
Via the Hungarian method this can be performed in $O(d^3)$ time; see \cite[\S17.3]{Schrijver03:CO_A}.
This is implemented in \polymake and can be used as shown in Listing~\ref{lst:tdet}.
\begin{lstlisting}[float, floatplacement=H,language=pmshell,caption={Computing a tropical determinant.},label=lst:tdet]
&\shelltrop& print tdet($A);
4
&\shelltrop& print tdet_and_perm($A);
4 <0 1 2>
&\shelltrop& print $A->elem(0,0) * $A->elem(1,1) * $A->elem(2,2);
4
\end{lstlisting}
The user can choose to only compute the value of $\tdet A$ or also one optimal permutation.
In the example from Listing~\ref{lst:tdet} that would be the identity permutation.

\section{Hypersurfaces}

A polyhedral complex is \emph{weighted}, if it is equipped with a function $\omega$, assigning integers to its maximal cells. A \emph{tropical cycle} is a weighted pure rational polyhedral complex $C$, such that each cell $C$ of codimension one satisfies a certain \emph{balancing condition}. The interested reader is referred to \cite{raumikhalkin, Tropical+Book}. Note that the latter only considers \emph{varieties}, which are cycles with strictly positive weights.
Tropical hypersurfaces and linear spaces are special cases of tropical varieties and our examples involve only these. 
We mention tropical cycles since the \polymake implementation is based on this concept.
Moreover, the tropical intersection theory, which we consider in Section \ref{section_intersection}, makes more sense in this general setting.

\subsection{Tropical hypersurfaces and dual subdivisions}\label{subsec:dual_subdivision}

Let
\[ F := \bigoplus_{a \in A} c_a \odot x^{\odot a} \in \T[x_0^\pm,\dots,x_n^\pm]\enspace,\]
be a \emph{tropical (Laurent) polynomial} with \emph{support} $A \subset \Z^{n+1}$, i.e., the coefficients $c_a$ are real numbers and $A$ is finite.
The \emph{tropical hypersurface} of $F$ is the set 
\[T(F) := \left\{ p \in \R^{n+1} \mid \text{the minimum in } F(p) \text{ is attained at least twice} \right\}\enspace.\]

Often we will assume that, for some $\delta\in\N$, we have $a_0+a_1+\dots+a_n=\delta$ for all $a\in A$.
This means that the tropical polynomial $F$ is \emph{homogeneous} (of degree $\delta$).
In this case for each point $p\in T(F)$ we have $p+\R\1\subseteq T(F)$.
Thus we usually consider the tropical hypersurface of a homogeneous polynomial as a subset of the quotient $\TT{n}$, which is called the \emph{tropical projective $n$-torus}. Note that one could also consider hypersurfaces in \emph{tropical projective space} $\T^n\backslash \{(\infty)^n)\} / \R\1$. However, from a computational point of view this incurs several challenges. In \polymake's implementation all tropical cycles live in the tropical projective torus and we will thus also adopt this viewpoint mathematically.

The \emph{dual subdivision} $\Delta(F)$ induced by $F$ is the collection of sets
\[\Delta_p := \left\{ \argmin_{a \in A} \{c_a \odot p^{\odot a}\}\right\} \subseteq A\enspace,\]
where $p$ ranges over all points in $\TT{n}$.
Instead of $\Delta(F)$ we also write $T(F)^*$.
Notice that, by definition $\Delta(F)$ is a set of subsets of $A$ whose union is $A$.
More precisely, $\Delta(F)$ is the combinatorial description of the regular subdivision of the support of $F$ induced by the coefficients.
We say that $T(F)$ is \emph{smooth} if the dual subdivision $\Delta(F)$ is unimodular, i.e., every maximal $\Delta_p$ is the vertex set of a unimodular simplex.

The tropical hypersurface of $F$ is the codimension one skeleton of the following subdivision of $\TT{n}$: We define two elements $p,p' \in \TT{n}$ to be equivalent, if $\Delta_p = \Delta_{p'}$. The equivalence classes are open polyhedral cones and their closures form a complete polyhedral complex, the \emph{normal complex} $\curly{D}(F)$.

As an example we consider the cubic polynomial
\begin{equation}\label{eq:nonvigeland_polynomial}
  \begin{split}
    F := &12 x_0^{\odot 3} \oplus (-131) x_0^{\odot 2} x_1 \oplus (-67) x_0^{\odot 2} x_2 \oplus (-9) x_0^{\odot 2} x_3 \oplus (-131) x_0 x_1^{\odot 2} \\ &\oplus (-129) x_0 x_1 x_2  \oplus (-131) x_0 x_1 x_3 \oplus (-116) x_0 x_2^{\odot 2} \oplus (-76) x_0 x_2 x_3\\ &\oplus (-24) x_0 x_3^{\odot 2} \oplus (-95) x_1^{\odot 3} \oplus (-108) x_1^{\odot 2} x_2 \oplus (-92) x_1^{\odot 2} x_3 \\ &\oplus (-115) x_1 x_2^{\odot 2} \oplus (-117) x_1 x_2 x_3 \oplus (-83) x_1 x_3^{\odot 2} \oplus (-119) x_2^{\odot 3} \\ &\oplus (-119) x_2^{\odot 2} x_3 \oplus (-82) x_2 x_3^{\odot 2} \oplus (-36) x_3^{\odot 3}
  \end{split}
\end{equation}
in four variables, i.e., the tropical hypersurface $V:=T(F)$ is a cubic surface in $\R^4/\R\1$.
This is constructed in Listing \ref{lst:cubic} along with the dual subdivision $V^*=\Delta(F)$.
\begin{lstlisting}[language=pmshell,caption={Computing a tropical cubic surface.},label=lst:cubic]
&\shelltrop& $F = toTropicalPolynomial("min(12+3*x0,-131+2*x0+x1,
  -67+2*x0+x2,-9+2*x0+x3,-131+x0+2*x1,-129+x0+x1+x2,
  -131+x0+x1+x3,-116+x0+2*x2,-76+x0+x2+x3,-24+x0+2*x3,-95+3*x1,
  -108+2*x1+x2,-92+2*x1+x3,-115+x1+2*x2,-117+x1+x2+x3,
  -83+x1+2*x3,-119+3*x2,-119+2*x2+x3,-82+x2+2*x3,-36+3*x3)");
&\shelltrop& $V = new Hypersurface<Min>(POLYNOMIAL=>$F);
&\shelltrop& print $V->DEGREE;
3
&\shelltrop& print $V->dual_subdivision()->N_MAXIMAL_CELLS;
27
\end{lstlisting}
The computation shows that $V$ is smooth.
Indeed, the support of $F$ are the lattice points of the scaled $3$-dimensional simplex $3\Delta_3$, whose normalized volume equals 27.
Since there are exactly 27 maximal cells, every single one must have volume 1 (and thus has to be a simplex as well).
We will come back to this example later to see that $V$ has some interesting enumerative properties.

While the entire design of the \polymake system follows the paradigm of object orientation there is a fundamental difference between an object of type \texttt{Matrix}, as in Listing~\ref{lst:matrix}, and an object of type \texttt{Hypersurface}, as in Listing~\ref{lst:cubic}.
Matrices form an example of a \emph{small object} class, while tropical hypersurfaces are \emph{big objects}.
To understand the difference it is important to know that, by design, \polymake employs both \texttt{Perl} and \texttt{C++} as main programming languages.
Essentially, the \texttt{C++} code deals with the computations for which speed matters.
The small objects belong to container classes which entirely live in the \texttt{C++} world.
On the \texttt{Perl} side this occurs as a mere reference, which is opaque.
Calling a member function on a small object from within the \polymake shell is always deferred to a corresponding \texttt{C++} function.
If new template instantiations occur for the first time this triggers just-in-time compilation.
The user experiences this as an occasional short time lag.
The newly compiled instantiation is kept in the \polymake folder in the user's home directory, such that it does not need to be compiled again.

Big objects are very different.
Technically, they entirely live in the \texttt{Perl} world.
More importantly, the user should think of them as technical realizations of actual mathematical objects, such as a tropical hypersurface.
For each big object class there is a certain number of \emph{properties} of which some subset is known at any given point in time.
In Listing~\ref{lst:cubic} the variable \texttt{\$V} is initialized as an object of type \texttt{Hypersurface<Min>} with the single property \texttt{POLYNOMIAL}, which is clearly enough to define a unique hypersurface.
The subsequent command prints a new property, the \texttt{DEGREE}, which is automatically derived from the input.
The essential idea is that this (and other properties computed on the way) will be kept and stored with the big object.
In this way it is avoided to repeat costly computations.
More details on \polymake's big object concept are found in \cite{Flexible+Object+Hierachies+in+polymake}.

\subsection{Product-Mix Auctions}

A fascinating connection between tropical geometry and economics was discovered by Baldwin and Klemperer \cite{Klemperer:2009, analysedemand}. They showed that the mechanics of \emph{product-mix auctions} can be modeled using tropical hypersurfaces. A more detailed analysis was given by Tran and Yu \cite{TranYu:1505.05737}, whose notation we mostly adopt. 

In a product-mix auction, several bidders (\enquote{agents}) compete for combinations of several goods, sold in discrete quantities. For example, one of the original motivations for this approach was, when the Bank of England wanted to auction off loans of funds during the financial crisis in 2007. These loans could be secured --- in various combinations --- against either weak or strong collateral. These two types of loans would be the goods in this case.

For our example, we will assume that there are two types of goods and only two agents.
Every agent now provides a \emph{utility function} $u^j: A^j \to \mathbb{R},\;j=1,2$, where $A^j \subseteq \mathbb{N}^2$ is the set of bundles of goods the agent is interested in.
Negative quantities could also be allowed, thus expressing an interest in selling the corresponding quantity.
The utility measures how valuable a bundle is to the agent.
Now, if the auctioneer fixes a price $p=(p_1,p_2)$, the agent will naturally be interested in the bundles which maximize her profit.
These bundles form the \emph{demand set}
\[D_{u^j}(p) := \argmax _{a \in A^j} \{u^j(a) - p \cdot a\}\enspace,\]
which depends not only on the price, but also on the choice of the utility function.
The \emph{aggregate demand} for the combined utilities $U=(u^1,u^2)$ is 
\[D_U(p) := \left\{  a^1 + a^2  \mid a^j \in D_{u^j}(p)\right\} \subseteq A^1+A^2\enspace.\]
Given an actual supply of $a \in \mathbb{N}^2$, the auctioneer will be interested in whether there exists a price such that all of the supply can be split between the agents such that every agent obtains a bundle which maximizes their profit, i.e., if there is a $p$ such that $a \in D_U(p)$.
In this case, we say that \emph{competitive equilibrium} exists at $a$.

In tropical language, every agent defines a hypersurface, corresponding to (the homogenization of) the tropical polynomial
\[F_j := \bigoplus_{a \in A^j} (- u^j(a)) \odot x^{\odot a}\enspace.\]
In this formulation, we see that $- F_j(p)$ is the maximal profit of agent $j$ at price $p$.
The tropical hypersurface $T(F_j)$ is the set of prices where the agent is indifferent between at least two bundles.

Let $f := f_1 \odot f_2$ be the product of the two polynomials and $A := A^1 + A^2$ its support.
Now, competitive equilibrium exists at a point $a \in A$, if and only if for some price vector $p$ the point $a$ is contained in the cell $\Delta_p$ of the dual subdivision $\Delta(F)$.

To illustrate, we compute Example 2 from \cite{TranYu:1505.05737}.
In Listing~\ref{lst:product} we define one hypersurface for each agent and a third one, $H$, which is the union.
\begin{lstlisting}[language=pmshell,caption={Constructing a tropical hypersurface from a product of polynomials},label=lst:product]
&\shelltrop& $H1 = new Hypersurface<Min>(
  MONOMIALS=>[[3,0,0],[2,0,1],[1,0,2],[0,1,2]], 
  COEFFICIENTS=>[0,-3,-5,-9]);
&\shelltrop& $H2 = new Hypersurface<Min>(
  MONOMIALS=>[[1,0,0],[0,1,0],[0,0,1]],
  COEFFICIENTS=>[0,-1,-1]);
&\shelltrop& $H = new Hypersurface<Min>(POLYNOMIAL=>
  $H1->POLYNOMIAL * $H2->POLYNOMIAL);
\end{lstlisting}
We need to homogenize the polynomial, so every bundle of goods has an additional coordinate in front.
Our goal is to determine the competitive equilibria.
Note that, since $H$ is a union of two tropical hypersurfaces, the dual subdivision $H^*$ is the common refinement of the dual subdivisions of the factors.
\begin{lstlisting}[language=pmshell,caption={The dual subdivision},label=lst:dual_subdivision]
&\shelltrop& $ds = $H->dual_subdivision();
&\shelltrop& $dehomog = $ds->POINTS->minor(All,~[0,1]);
&\shelltrop& $cells = transpose($ds->MAXIMAL_CELLS);
\end{lstlisting}
The monomials of the tropical polynomial defining $H$ arise as the vertices of the cells of $H^*$.
To reinterpret them as bundles we dehomogenize, i.e., we strip the first two coordinates; the first one equals one (since ordinary points are homogenized), while the second one equals zero (due to the tropical homogenization).
For each bundle or monomial we can now print the number of cells containing it.
\begin{lstlisting}[language=pmshell,caption={Checking all bundles},label=lst:competitive]
&\shelltrop& for (my $i=0; $i<$ds->N_POINTS; ++$i) {
  print $dehomog->row($i), ": ", $cells->row($i)->size(), "\n" }
0 0: 2
1 0: 1
1 3: 2
0 1: 2
1 1: 0
0 2: 2
1 2: 5
2 2: 2
0 3: 1
\end{lstlisting}
Indeed, we see that every bundle, except for $(1,1)$, is in at least one cell of the dual subdivision.
That is, competitive equilibrium exists precisely at the nine remaining bundles.

\section{Moduli of Tropical Plane Curves}

So far we investigated individual examples of tropical hypersurfaces.
Now we will look into families which are obtained by varying the coefficients.
To this end we start out with a point configuration and a given subdivision $\Sigma$.
The goal is to determine all possible tropical polynomials $F$ with $\Delta(F)=\Sigma$.
Our example computation will deal with a planar point configuration, and hence the tropical hypersurfaces $T(F)$ will be tropical plane curves.
These objects stood at the cradle of tropical geometry; see, in particular, Mikhalkin~\cite{mikhalkinenumerative}.
Yet the study of their moduli spaces is more recent \cite{BJMS:2015}, and this is the direction where we are heading here; see also \cite{BirkmeyerGathmann:1412.3035}.

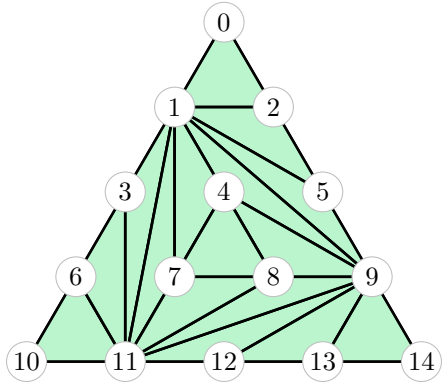
\begin{figure}[bh]
  \centering
  \input{smoothKleinTriangulation.tikz}
  \caption{Unimodular triangulation of $4\Delta_2$}
  \label{fig:klein:triangulation}
\end{figure}

The Listing~\ref{lst:klein:triangulation} shows \polymake code to visualize a triangulation of 15 points in the affine hyperplane $\sum x_i = 4$ in $\R^3$; see Figure~\ref{fig:klein:triangulation}.
For technical reasons the set of points is converted into a matrix with leading ones.
We start out by switching the application.
\begin{lstlisting}[float,floatplacement=H,language=pmshell,caption={Constructing and visualizing a triangulation.},label=lst:klein:triangulation]
&\shelltrop& application "fan";
&\shellfan& $points = [ [0,0,4],[1,0,3],[0,1,3],[2,0,2],[1,1,2],
 [0,2,2], [3,0,1],[2,1,1],[1,2,1], [0,3,1],[4,0,0],[3,1,0],
 [2,2,0],[1,3,0],[0,4,0] ];
&\shellfan& $triangulation = [[0,1,2],[9,11,12],[9,12,13],
 [9,13,14],[1,2,5],[6,10,11],[3,6,11],[1,5,9],[1,3,11],
 [8,9,11],[1,4,9],[1,7,11],[7,8,11],[4,8,9],[1,4,7],[4,7,8]];
&\shellfan& $pointMatrix =
 (ones_vector<Rational>(15)) | (new Matrix<Rational>($points));
&\shellfan& $Sigma = new SubdivisionOfPoints(
 POINTS=>$pointMatrix, MAXIMAL_CELLS=>$triangulation);
&\shellfan& $Sigma->VISUAL;
\end{lstlisting}

Let $P$ be a finite set of points in $\R^d$, and let $\Sigma$ be a \emph{(polytopal) subdivision} of $P$, i.e., $\Sigma$ is a polytopal complex whose vertices form a subset of $P$ and which covers the convex hull of $P$.
The \emph{secondary cone} $\seccone \Sigma$ is the topological closure of the set of lifting functions on the set $P$ which induce $\Sigma$.
Scaling any lifting function in $\seccone \Sigma$ by a positive number does not lead outside, and neither does adding two such lifting functions.
We infer that $\seccone \Sigma$, indeed, is a cone.
Each cell of codimension one in $\Sigma$ yields one linear inequality, and these give an exterior description of the secondary cone.
This means that the secondary cone is polyhedral.
It is of interest to determine the rays of $\seccone \Sigma$.
For our example this is accomplished in Listing \ref{lst:klein:seccone}.

\begin{lstlisting}[language=pmshell,caption={Analyzing the secondary cone.},label=lst:klein:seccone]
&\shellfan& $sc=$Sigma->secondary_cone();
&\shellfan& print $sc->N_RAYS;
12
&\shellfan& for (my $i=0; $i<$sc->N_RAYS; ++$i) {
       my $c=new SubdivisionOfPoints(
              POINTS=>$pointMatrix,WEIGHTS=>$sc->RAYS->[$i]);
       print $i, ":", $c->N_MAXIMAL_CELLS, " ";
      }
0:2 1:2 2:2 3:2 4:2 5:2 6:2 7:2 8:2 9:3 10:3 11:3
\end{lstlisting}

\begin{figure}[bh]
  \centering
  \begin{minipage}{.23\linewidth}
    \input{smoothKleinTriangulation-ray0.tikz}
  \end{minipage}
  \hfill
  \begin{minipage}{.23\linewidth}
    \input{smoothKleinTriangulation-ray2.tikz}
  \end{minipage}
  \hfill
  \begin{minipage}{.23\linewidth}
    \input{smoothKleinTriangulation-ray8.tikz}
  \end{minipage}
  \hfill
  \begin{minipage}{.23\linewidth}
    \input{smoothKleinTriangulation-ray11.tikz}
  \end{minipage}
  \caption{The coarsest subdivisions of the rays labeled 0, 2, 8 and 11.}
  \label{fig:klein:rays}
\end{figure}
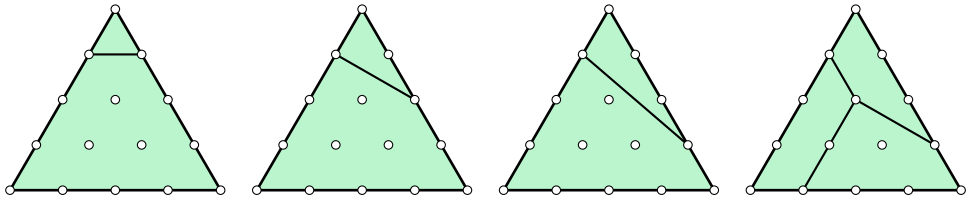

The rays of $\seccone \Sigma$ induce those \emph{coarsest subdivisions} of the point set $P$ from which $\Sigma$ arises as their common refinement.
In our example the secondary cone has 12 rays, which come in no particular order.
In Listing \ref{lst:klein:seccone} we list the index of each ray (from 0 to 11) with the number of maximal cells in the corresponding coarsest subdivision.
Throughout these numbers are either two or three, from which one can tell right away that the former are $2$-splits, while the latter are $3$-splits; see Herrmann~\cite{Herrmann:2011}.
The 12 rays come in four orbits, with respect to the symmetry group of $P$ which fixes $\Sigma$.
The order of that group is three.
There are three orbits of $2$-splits, represented by 0, 2 and 8, and one orbit of $3$-splits, represented by 11.
These are shown in Figure~\ref{fig:klein:rays}.

One specific lifting function on the 15 points which yields our example triangulation, $\Sigma$, is shown in Listing~\ref{lst:klein:curve}.
We use it as the coefficient vector of a tropical polynomial, and this defines a tropical hypersurface, which we call $C$,
To verify that this vector, indeed, lies in the relative interior of the secondary cone of $\Sigma$ we can compute the scalar products with all facet normal vectors.
\begin{lstlisting}[language=pmshell,caption={A tropical plane curve.},label=lst:klein:curve]
&\shellfan& application "tropical";
&\shelltrop& $C = new Hypersurface<Min>(MONOMIALS=>$points,
 COEFFICIENTS=>[6, 0,3, 1,-1/3,1, 3,-1/3,-1/3,0, 6,0,1,3,6]);
&\shelltrop& $ratCoeff = new Vector<Rational>($C->COEFFICIENTS);
&\shelltrop& print $sc->FACETS * $ratCoeff;
4 4 8/3 4 4 4 4 8/3 8/3 4/3 4/3 4/3
&\shelltrop& print $sc->LINEAR_SPAN * $ratCoeff;
\end{lstlisting}
The fact that all these numbers are positive serves as a certificate for strict containment; the actual values do not matter.
For a general subdivision, which is not a triangulation, there are additional linear equations to be checked which describe the linear span.
There are no such equations in this case, which is why the last command has no output.
Note that prior to computing the scalar products it is necessary to explicitly convert the lifting function into a vector with rational coefficients.
This is unavoidable since we want to use the ordinary scalar multiplication here, not the tropical one.

\begin{figure}[ht]
  \centering

  \input{smoothKleinCurve.tikz}

  \caption{Tropical plane quartic of genus three.}
  \label{fig:klein:curve}
\end{figure}
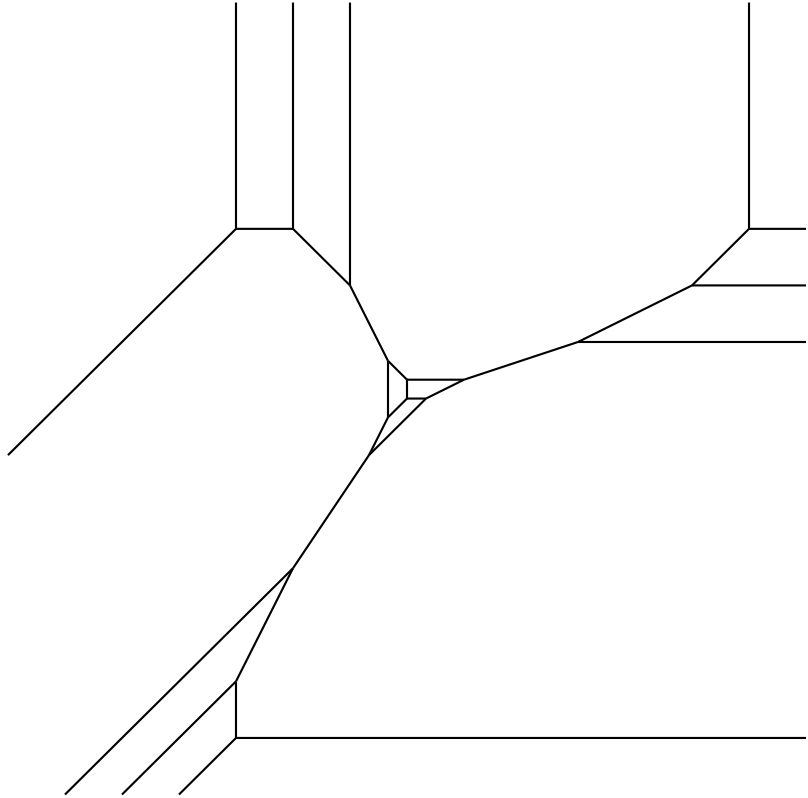

\begin{figure}[ht]
  \begin{minipage}[b]{.47\linewidth}\raggedright
    \input{smoothKleinSkeleton-A.tikz}
  \end{minipage}
  \hfill
  \begin{minipage}[b]{.47\linewidth}\raggedleft
    \input{smoothKleinSkeleton-B.tikz}
  \end{minipage}

  \caption{Essential part of a tropical plane quartic (with edge lengths) and skeleton (with edge labels).}
  \label{fig:klein:skeleton}
\end{figure}
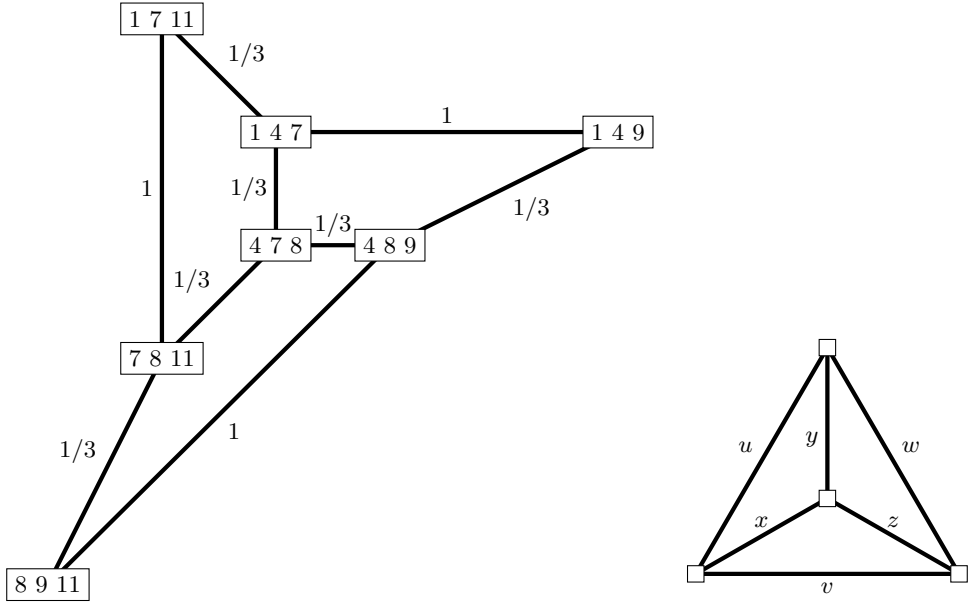

The curve $C$ defined in Listing~\ref{lst:klein:curve} and shown in Figure~\ref{fig:klein:curve} is a tropical plane quartic of genus three:
the Newton polytope is the dilated standard simplex $4\Delta_2$, and this has exactly three interior lattice points, all of which are used in the triangulation $\Sigma$. 
\begin{lstlisting}[language=pmshell,caption={Degree and genus of a tropical plane curve.},label=lst:klein:degree]
&\shelltrop& print $C->DEGREE;
4
&\shelltrop& print $C->GENUS;
3
\end{lstlisting}
In the sequel we want to locate $C$ in the moduli space of tropical plane curves of genus three.
This will go hand in hand with our initial goal to determine all such curves which fit the initial triangulation $\Sigma$ of $4\Delta_2$.
To this end we determine the (lattice) length of each edge of $C$, considered as a one-dimensional ordinary polytopal complex.
The integer lattice $\Z^3$ induces sublattices on the line spanned by any edge of $\Sigma$.
Since each edge of $C$ is dual to an edge of $\Sigma$ we can measure its length with respect to that sublattice.
The result is shown in Listing~\ref{lst:klein:moduli}.
The 30 edges are labeled from 0 through 29; again they do not come in any particular order.
\begin{lstlisting}[float, floatplacement=H,language=pmshell,caption={Moduli of a tropical plane curve.},label=lst:klein:moduli]
&\shelltrop& print labeled($C->CURVE_EDGE_LENGTHS);
0:inf 1:inf 2:1 3:inf 4:1 5:1 6:inf 7:inf 8:1/3 9:1/3 10:1/3
 11:1/3 12:2/3 13:1/3 14:inf 15:1 16:inf 17:2/3 18:1 19:1 20:1
 21:1/3 22:2/3 23:inf 24:1 25:inf 26:1 27:inf 28:inf 29:inf
&\shelltrop& $C->VISUAL(LengthLabels=>"show");
\end{lstlisting}
As a generic curve of degree four, the curve $C$ has four edges of infinite length in each of the three coordinate directions.
Contracting these infinite edges including the edges of finite length which ``lead'' to them yields the ``essential part'' of $C$.
A picture of this with the remaining edge lengths is given in Figure~\ref{fig:klein:skeleton} (left), the remaining vertices are labeled with the corresponding triangles of $\Sigma$.
In the essential part of $C$ we have vertices of degree two or three.
Joining edges by omitting those of degree two gives the \emph{combinatorial skeleton} of $C$, which is a planar graph with $2g - 2$ vertices and $3g - 3$ edges, where $g$ is the genus, i.e., $g=3$ in our case.
The name comes about from its (loose) connection to the Berkovich skeleton of the analytification of a smooth complete curve; see \cite{GublerRabinoffWerner:1508.01179}.
The joined edges receive the sum of the lengths of the original edges of the curve, and this way we arrive at a metric graph.
The lengths of the skeleton edges are the \emph{moduli} of $C$.
In this case the skeleton is the \emph{honeycomb graph} of genus three, denoted as ``$(000)$'' in~\cite{BJMS:2015}.
Figure~\ref{fig:klein:skeleton} (right) shows the skeleton of $C$ with the edge labels as in \cite[Fig.~4]{BJMS:2015}; the moduli are
\[
u \ = \ v \ = \ w \ = \ 1 + \frac{1}{3} \ = \ \frac{4}{3} \quad \text{and} \quad x \ = \ y \ = \ z \ = \ \frac{1}{3} \enspace ,
\]
and this agrees with \cite[Thm.~5.1]{BJMS:2015}.


Without going through all the remaining computations explicitly we now want to sketch how the secondary cone of $\Sigma$ contributes to the moduli space of tropical plane curves of genus three.
Via measuring lattice lengths of edges and properly attributing their contributions to the $3g-3$ edges of a fixed skeleton curve, $\seccone \Sigma$ is linearly mapped to a cone in $\R^{3g-3}$; see \cite[pp.~3ff]{BJMS:2015} for the details.
We apply this procedure to the 12 rays computed in Listing~\ref{lst:klein:seccone} and visualized in Figure~\ref{fig:klein:rays}.
Comparing the dual pair of pictures in Figures~\ref{fig:klein:triangulation} and~\ref{fig:klein:curve} we can see that the three orbits of $2$-splits do not contribute anything to essential part of the curve.
This means they are mapped to zero in the moduli cone.
On the other hand, e.g., the $3$-split corresponding to ray~11 corresponds to a curve with moduli $v=1$, $x=z=\frac{1}{3}$ and $u=w=y=0$.
The three rays in this orbit (with labels 9, 10 and 11) span a $3$-dimensional cone in the moduli space, and this is the image of $\seccone \Sigma$ under the linear map described above.
The coefficient vector which defines the curve is, up to scaling, the sum of those three rays.
That is, it sits right in the center of that moduli cone.

The moduli space of tropical plane curves of genus $3$ is generated from the secondary cones of all 1278 unimodular triangulations (up to symmetry) of the dilated triangle $4\Delta_2$, plus an additional contribution from hyperelliptic curves.
The dimension of the entire moduli space is six. Each of the 1278 triangulations contributes a moduli cone whose dimension lies between three and six; see \cite[Table~1]{BJMS:2015}.
For higher genus it is necessary to consider several polygons; see \cite{CastryckVoight:2009}.
Since that moduli space describes \emph{isomorphism classes} of tropical curves it is necessary to take symmetries (of the triangulations and the skeleta) into account.
This entails that the global structure is \emph{not} a polyhedral fan but rather a quotient structure called \emph{stacky fan}; see \cite{BrannettiMeloViviani:2011}.




\section{Grassmannians,  Linear Spaces and Matroids}

There are many equivalent ways to define matroids, and we recommend that the interested reader look at \cite{Oxley:2011,White:1986}. For our purposes the criterion by Gel$'$fand et.\ al.\cite{GelfandEtAl:1987} is the most convenient: Let $M \subseteq \binom{[n]}{k}$ be a set of subsets of $[n]=\{1,2,\dots,n\}$ of size $k$. To this collection we can associate a polytope
$$P_M := \conv\{v_B := \sum_{i \in B} e_i \mid B \in M\}\enspace,$$
where $e_i$ is a standard basis vector.
If $M = \binom{[n]}{k}$ consists of all $k$-sets, we call $P_M =: \Delta(k,n)$ a \emph{hypersimplex}. We say that $M$ is a \emph{matroid} of rank $k$ on $[n]$, if every edge of $P_M$ is parallel to $e_i - e_j$ for some $i \neq j$. We call $P_M$ the \emph{matroid (basis) polytope} of $M$ and the elements of $M$ its \emph{bases}. We say that $M$ is \emph{loopfree}, if every element of $[n]$ is contained in some basis.

\subsection{The tropical Grassmannian}

The \emph{tropical Grassmannian} $\TGr(k,n)$ was first studied in detail by Speyer and Sturmfels \cite{SpeyerSturmfels04}. It is the tropicalization of the (complex) Grassmannian $\Gr(k,n)$, intersected with the torus. It also parametrizes tropicalizations of uniform linear spaces. This can be seen in the following fashion: Every element of $\TGr(k,n)$ is a \emph{tropical Plücker vector} $p \in \mathbb{R}^{\binom{n}{k}}$. Equivalently, we can view it as a height function on the hypersimplex $\Delta(k,n) \subseteq \mathbb{R}^n$.
The set of all tropical Pl\"ucker vectors is the \emph{Dressian} $\Dr(k,n)$.  In general, $\TGr(k,n)$ is a proper subset of $\Dr(k,n)$.  Throughout the following we assume that $k\leq n$.

This height function thus induces a regular subdivision of $\Delta(k,n)$. The fact that $p$ is a Plücker vector implies (but is generally not equivalent) that this subdivision is matroidal, i.e.\ every cell is again a matroid basis polytope. That is, if the set $\{v_{B_1},\dots,v_{B_k}\}$ comprises the vertices of a cell, then $\{B_1,\dots,B_k\}$ is the set of bases of a matroid.

The combinatorics of $\TGr(3,6)$ were studied in detail in \cite{drawtropicalplanes}. The authors compute that there are seven combinatorial types of generic uniform linear spaces, basically encoded in their bounded complexes. 
As an example, we want to consider a particular vector $p \in\mathbb{R}^{\binom{6}{3}}$.
We start out with analyzing the combinatorics, and we turn to algebraic computations later.

\begin{lstlisting}[language=pmshell,caption={Computing the tight span of a tropical Plücker vector},label=lstpluecker]
&\shelltrop& $Delta=polytope::hypersimplex(3,6);
&\shelltrop& $p=new Vector<Int>(
  [0,0,3,1,2,1,0,1,0,2,2,0,3,0,4,1,2,2,0,0]);
&\shelltrop& $tlinear=new fan::SubdivisionOfPoints(
  POINTS=>$Delta->VERTICES, WEIGHTS=>$p);
&\shelltrop& print $tlinear->TIGHT_SPAN->MAXIMAL_POLYTOPES;
{0 4}
{1 5}
{1 2 3 4}
&\shelltrop& $bases = new IncidenceMatrix(
  matroid::uniform_matroid(3,6)->BASES);
&\shelltrop& @subdiv_bases = map { 
  new Array<Set>(rows($bases->minor($_,All))) 
  } @{$tlinear->MAXIMAL_CELLS};
&\shelltrop& print join(",", map { 
  matroid::check_basis_exchange_axiom($_) } @subdiv_bases );
1,1,1,1,1,1
\end{lstlisting}

Now we want to verify that the vector $p$, indeed, lies in $\TGr(3,6)$.  It is known that for $(k,n)=(3,6)$ the tropical Grassmannian and the Dressian agree as sets.  However, in general, this does not hold.
We employ \polymake 's interface to \Singular, and we switch to the application \texttt{ideal}.

\begin{lstlisting}[float,floatplacement=H,language=pmshell,caption={Computing a generalized initial ideal of the Plücker ideal.},label=lstinitial]
&\shelltrop& application "ideal";
&\shellideal& $I=pluecker_ideal(3,6);
&\shellideal& $pp=new Vector<Int>(5*ones_vector(20)-$p);
&\shellideal& $J=new Ideal(GENERATORS=>
  $I->GROEBNER(ORDER_VECTOR=>$pp)->INITIAL_FORMS);
&\shellideal& print $J->contains_monomial();
0
\end{lstlisting}
A few explanations are in order.  The \emph{Pl\"ucker ideal} $I(k,n)$ describes the algebraic relations among the $k{\times}k$-minors of a general $k{\times}n$-matrix.  This is an ideal in the polynomial ring over the integers with $\tbinom{n}{k}$ indeterminates, each of which encodes a choice of $k$ columns to specify one such minor.  There is a purely combinatorial description  of the reverse lex Gr\"obner basis of $I(k,n)$, and this is what is computed by \polymake directly; see \cite[Chapter~3]{Sturmfels:2008}. This function is also implemented in \Macaulay \cite{M2}.  The tropical variety $\TGr(k,n)$ arises as a subfan of the Gr\"obner fan of $I(k,n)$.  Yet it is common that the interpretation of the vectors in the Gr\"obner fan refer to maximization, while our choice for regular subdivisions relies on minimization.  This entails that we need to swap from the tropical Pl\"ucker vector $p$ to its negative. For technical reasons \Singular requires such weight vectors to be positive, which is why we consider $p'=5\cdot\1-p$; notice that $\tbinom{6}{3}=20$.  The condition for $p'$ to lie in $\TGr(3,6)$ is that the ideal which is generated by the leading forms of $I(3,6)$ with respect to $p'$ does not contain any monomial.
\begin{remark}
  Our choice for $p$ or rather $p'$ corresponds to a generic tropical $2$-plane in $5$-space of type EEFG in the notation of \cite{SpeyerSturmfels04}.
\end{remark}

\subsection{Tropical linear spaces}

So far, we have only considered special Pl\"ucker vectors --- in the sense that the underlying matroid is a \emph{uniform} matroid. There is a general theory of valuated matroids, originally developed by Dress and Wenzel \cite{dresswenzel}. 

\begin{definition}
 A \emph{valuated matroid} $(M,v)$ is a matroid $M$ together with a function $v$ from the set of its bases to the real numbers such that the induced regular subdivision on the matroid basis polytope of $M$ is matroidal, i.e.\ every cell is a matroid polytope again.
 
 We denote the regular subdivision by $\Delta(M,v)$. It has a normal complex $\curly{D}(M,v)$; see Section \ref{subsec:dual_subdivision}. The \emph{tropical linear space} $B(M,v)$ associated to $(M,v)$ is the subcomplex of $\curly{D}(M,v)$ consisting of all faces whose corresponding dual cell if the basis polytope of a loopfree matroid.
\end{definition}

If the matroid is the uniform matroid $U_{k,n}$, then $v$ is an element in the Dressian $\Dr(k,n)$. If, moreover, it is realizable, then $v$ lies in the tropical Grassmannian $\TGr(k,n)$. 
Tropical linear spaces play an important role in tropical geometry. They are exactly the tropical varieties of degree 1 \cite{finkchow}. If the valuation $v \equiv 0$ is trivial, the tropical linear space is a fan, also called the \emph{Bergman fan} $B(M)$. These fans are the basic building blocks for \emph{smooth} tropical varieties; see the discussion in \ref{section_intersection}. Also, the tropical homology of a Bergman fan encodes the Orlik-Solomon algebra of the matroid \cite{zharkov}.

The definition itself already suggests an algorithm to compute a tropical linear space: Find all cells of the matroidal subdivision corresponding to a loopfree matroid and compute its cell in the normal complex. For this, \polymake uses a variant of Ganter's algorithm\cite{Ganter:1984,Ganterreuter:1991}, which computes, in fact, the full face lattice of the tropical linear space.

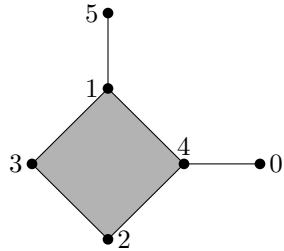
\begin{figure}[ht]
\centering
 \begin{tikzpicture}
  \fill[light_gray, draw = black,text=black] (-1,0) node[left]{3} -- (0,1) node[left] {1} -- (1,0) node[above]{4} -- (0,-1) node[right]{2} -- cycle;
  \draw (0,1) -- (0,2) node[left]{5};
  \draw (1,0) -- (2,0) node[right]{0};
  \foreach \x in {(-1,0),(0,1),(1,0),(0,-1),(0,2),(2,0)}
      {\fill[black] \x circle (2pt);}
 \end{tikzpicture}
 \caption{The combinatorial type of the tropical linear space computed in Listing \ref{lstlinearspace}.}\label{fig:linearspace}
\end{figure}

As an example, we want to compute the tropical linear space associated to the (uniform) Pl\"ucker vector we considered in Listing \ref{lstpluecker}. We then compute the complex of its bounded faces to confirm that its combinatorics matches the one in Figure \ref{fig:linearspace}. In this example, the bounded faces are identified as those, which do not contain any rays. The latter are stored in the property \texttt{FAR\_VERTICES}.

\begin{lstlisting}[language=pmshell,caption={Computing a tropical linear space.},label=lstlinearspace]
&\shellideal& application "tropical";
&\shelltrop& $mat = new matroid::ValuatedMatroid<Min>(
  BASES=>\@{rows($bases)}, VALUATION_ON_BASES=>$tlinear->WEIGHTS, 
  N_ELEMENTS=>6);
&\shelltrop& $tl= linear_space($mat);
&\shelltrop& @bounded = map { 
  grep { ($_ * $tl->FAR_VERTICES)->size == 0} @{rows($_)} 
  } @{$tl->CONES};
&\shelltrop& print join(",",@bounded);
{0},{1},{2},{3},{4},{5},{0 4},{1 4},{2 4},{1 3},{2 3},{1 5},
  {1 2 3 4}
\end{lstlisting}

In the case of Bergman fans there are two more methods to compute the space:

\subsubsection*{Cyclic fan structure}

Rinc\'{o}n studies the \emph{cyclic Bergman fan} of a matroid \cite{rinconcompute}, which is a refinement of the polyhedral structure we have considered so far. It relies heavily on computing \emph{fundamental circuits} and is particularly fast in the case of matrix matroids, where these computations can be carried out by standard linear algebra methods.

\subsubsection*{Order complex of the lattice of flats} Ardila and Klivans \cite{ArdilaKlivans:2006} proved that the order complex of the lattice of flats of a matroid $M$ can be realized as a fan which is supported on $B(M)$. That is, we obtain a refinement of $B(M)$, such that every ray corresponds to a flat and every cone to a chain of flats. Due to the potentially large number of flats, this is naturally not a very efficient method to compute a tropical linear space. Nevertheless, this is still feasible for small matroids and this particular subdivision is often useful.

In Listing \ref{lstbergman}, we compute both these fans for the Fano matroid, plus the linear space of its trivial valuation. We see that in this particular case, all these fans are actually the same.
\begin{lstlisting}[language=pmshell,caption={Computing a Bergman fan in various ways.},label=lstbergman]
&\shelltrop& $fano = matroid::fano_matroid();
&\shelltrop& $cyclic = matroid_fan<Min>($fano);
&\shelltrop& $linear = linear_space(
  matroid::trivial_valuation<Min>($fano));
&\shelltrop& $order = matroid_fan_from_flats<Min>($fano);
&\shelltrop& print join(",", map {
  $_->N_MAXIMAL_POLYTOPES } ($cyclic,$linear,$order));
21,21,21
\end{lstlisting}

\section{Intersection Theory} \label{section_intersection}

The basics for a tropical intersection theory were already laid out in \cite{mikhalkinenumerative}. They are closely related to the \emph{fan displacement rule} by Fulton and Sturmfels \cite{fultonsturmfels}. The upshot is that, to intersect two tropical cycles in $\TT{n}$, one of them is shifted in a generic direction until they intersect transversely. The actual intersection product is the limit of the transversal product when shifting back. 

This definition is of course hardly suitable for computations. Instead, \atint uses the criterion by Jensen and Yu \cite{jensenyu} for computations. This is a local criterion, which states that a point in the set-theoretic intersection of two cycles $X$ and $Y$ is also in the intersection product, if and only if the Minkowski sum of the local fans at this point is full-dimensional (a more detailed account of various definitions of the intersection product in the tropical torus can be found in \cite{hatint}).

The situation becomes more difficult when the ambient variety is not the whole projective torus. As in the algebraic case, even the theory is only fully understood in the case of \emph{smooth varieties}. This is a tropical cycle with positive weights, which is everywhere locally isomorphic to the Bergman fan of a matroid (where an isomorphism of fans is a map in $\textnormal{GL}_n(\Z)$ respecting weights). Every tropical linear space is smooth by definition. For hypersurfaces, this definition of smoothness coincides with the one given in \ref{subsec:dual_subdivision}.

Since intersection products should be local, it is enough to understand how to compute them in Bergman fans of matroids. In this case, there are two equivalent definitions by Shaw \cite{shawmatroidal} and Fran\c{c}ois and Rau \cite{francoisrau}. In general, both are not very computation-friendly. The first is a recursive procedure using projections and pull-backs. The second employs the idea of \enquote{cutting out the diagonal}, i.e.\ writing down rational functions whose consecutive application gives the diagonal of $B(M) \times B(M)$. Due to the large dimension and the quadratic increase in the number of cones, this method does not perform well. 

However, Shaw's method does simplify nicely in the case of surfaces (i.e.\ matroids of rank 3), where the intersection product can be computed in a nice combinatorial manner \cite[Section 4]{shawmatroidal}. Using a smoothness detection algorithm, which was implemented in \atint by Dennis Diefenbach, this enabled us to write a procedure which computes intersection products of cycles in smooth surfaces. The upshot of the detection algorithm is that in the case of surfaces, one can basically do a brute force search over all possible ways to assign flats to the rays of the fan.

As a demonstration of both intersection in the tropical torus and in a smooth surface, we will consider the hypersurface defined by the polynomial (\ref{eq:nonvigeland_polynomial}) in \ref{subsec:dual_subdivision}. The reader will have to believe (or verify) that this hypersurface contains the standard tropical line with apex $(0,0,0,0)$ or, in other words, the Bergman fan $B$ of $U_{2,4}$. We want to compute the self-intersection of $B$ in $V$. Also, we will calculate the 3-fold intersection of $V$ in the torus. We already know from the Tropical Bernstein Theorem \cite[Theorem 4.6.8]{Tropical+Book} that this is the lattice volume of $3 \Delta_3$.

\begin{lstlisting}[language=pmshell,caption={Self-intersection in a smooth surface},label=lstdegenerate]
&\shelltrop& print intersect(intersect($V,$V),$V)->DEGREE;
27
&\shelltrop& $B = matroid_fan<Min>(matroid::uniform_matroid(2,4));
&\shelltrop& print intersect_in_smooth_surface($V,$B,$B)->DEGREE;
-1
\end{lstlisting}

This seems to tie in nicely with the classic fact that a line in a smooth cubic has self-intersection $-1$. Hence we want to take a closer look at that situation.

\subsection{Lines in tropical cubics}

In algebraic geometry it is a well-known fact that any smooth cubic surface in $\mathbb{P}^3$ contains exactly 27 lines.
It is known that the incidence structure arising from the 27 lines and their 45 points of intersection is the unique generalized quadrangle of order $(4,2)$; see \cite[\S3]{ThasVanMaldeghem:2006} for details and related constructions.
For instance, it is known that any line intersects exactly ten other lines, and for any two disjoint lines there are five lines that intersect both of them.
Furthermore, as mentioned before, they all have self-intersection $-1$.

In tropical geometry, the situation is much more complicated --- or, possibly, much more interesting, depending on your point of view. First of all, we need to establish what a tropical line is:

\begin{definition}
 A \emph{tropical line} in $\TT{n}$ is the tropical linear space of a valuation on $U_{2,n+1}$.
\end{definition}

Now the peculiarities begin with the fact that a smooth tropical cubic surface may actually contain \emph{families} of tropical lines. A first systematic study of this problem was undertaken by Vigeland in \cite{Vigeland}. He provided an example of a secondary cone of $3\Delta_3$, such that any general element of that cone defines a tropical cubic which contains exactly 27 lines. He also gave a list of possible \emph{combinatorial types}, which describe how a line can lie in a tropical surface.
Our example, which does not occur in Vigeland's list, was found via a systematic search through the secondary fan of $3\Delta_3$.
Here we just show our example, while its complete analysis is beyond the scope of the present paper,

Other approaches focused on counting the lines \enquote{in the correct manner}, e.g., by showing that certain lines could not be \emph{realized} as tropicalizations of lines in a cubic surface \cite{BrugalleShaw, BogartKatz}. The paper \cite{RenShawSturmfels} finds the 27 lines as trees in the boundary of the tropicalization.
The interest in this particular problem stems from the fact that it provides a nontrivial, yet computationally feasible testing ground for studying the problem of (relative) \emph{realizability}. It may also serve as an indicator of what possible additional structure should be associated to a tropical variety, i.e.\ in bold terms, what a \emph{tropical scheme theory} could be; see for example \cite{giansiracusas, maclaganrincon}. 

\subsection{Vigeland's missing type}

We will reconsider the polynomial (\ref{eq:nonvigeland_polynomial}) from \ref{subsec:dual_subdivision}. 
We compute the list of lines (and families thereof) in the tropical hypersurface corresponding to $f$:

\begin{lstlisting}[language=pmshell,caption={Computing lines in a tropical cubic surface.},label=lstlines]
&\shelltrop& $L = lines_in_cubic($F);
[Output omitted]
&\shelltrop& print $L->N_ISOLATED,", ",$L->N_FAMILIES,"\n";
26, 3
\end{lstlisting}

This demonstrates that there are in fact 26 isolated lines and three different families of lines in the tropical hypersurface $V$ defined by $f$.
It can be shown that small random changes to the coefficients do not affect the combinatorics of the lines in the corresponding cubic.
In particular, this contradicts Conjecture 1 in \cite{Vigeland} that a general cubic contains exactly 27 lines.
In Vigeland's terminology, a general cubic with a fixed dual subdivision corresponds to a dense open subset in the euclidean topology in the secondary cone.
A full discussion and formal proof of our claims will be deferred to another paper.

Now we want to take a closer look at the families. One can ask \polymake for a picture of the families using

\begin{lstlisting}[language=pmshell,caption={Visualizing families in  the cubic.},label=lstvisual]
&\shelltrop& @rp = map { $_->representative() } $L->all_families;
&\shelltrop& visualize_in_surface($V, @rp);
\end{lstlisting}

This produces a (more or less generic) representative of each family and visualizes them together in the hypersurface.

\begin{figure}
\begin{tabular}{c c}
 \includegraphics[scale=.2]{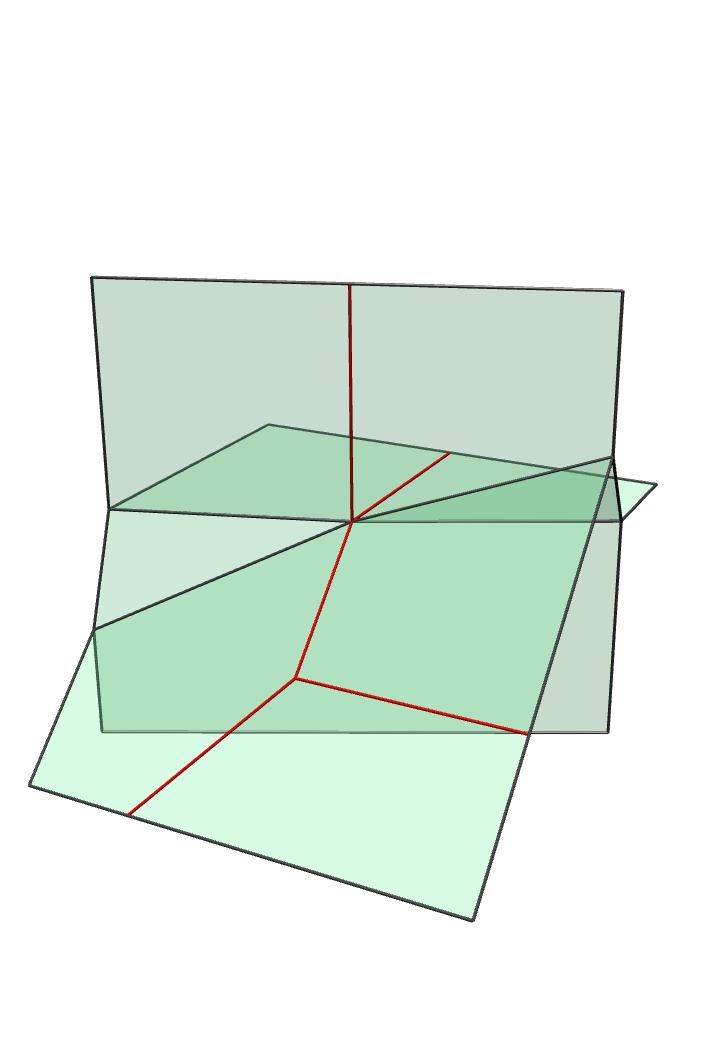} & 
 \includegraphics[scale=.2]{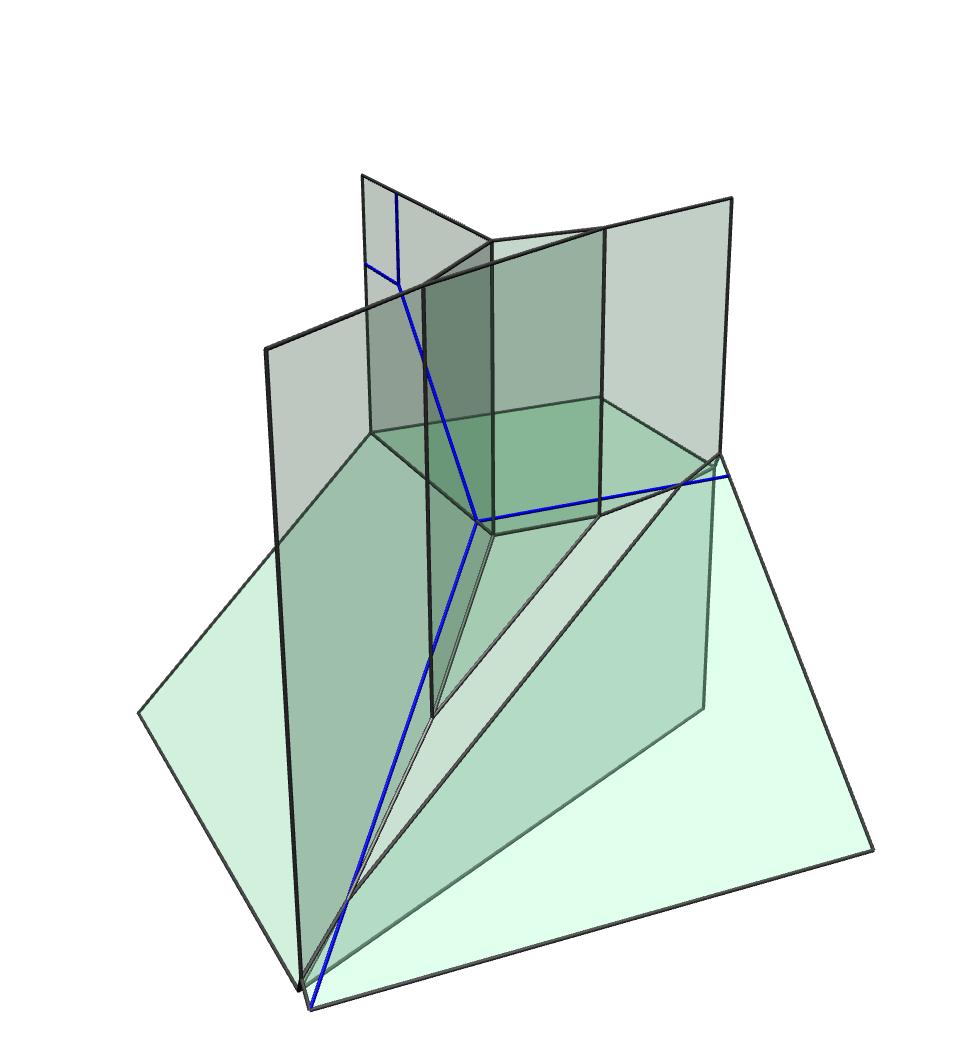} \\
 \begin{tikzpicture}[scale=.3]
  \draw (-1,1) -- (0,0) -- (-1,-1);
  \draw (0,0) -- (1,0);
  \draw (2,1) -- (1,0) -- (2,-1);
  \fill[black] (0,0) circle (5pt);
 \end{tikzpicture} &
 \begin{tikzpicture}[scale=.3]
  \draw (-1,1) -- (0,0) -- (-1,-1);
  \draw (0,0) -- (1,0);
  \draw (2,1) -- (1,0) -- (2,-1);
  \fill[black] (1.5,.5) circle (5pt);
  \fill[black] (1.5,-.5) circle (5pt);
  \draw[line width=2pt] (0,.5) -- (0,-.5);
 \end{tikzpicture}
\end{tabular}
\caption{Local pictures of representatives of families in the cubic and their combinatorial type in Vigeland's notation.}\label{fig_families}
\end{figure}
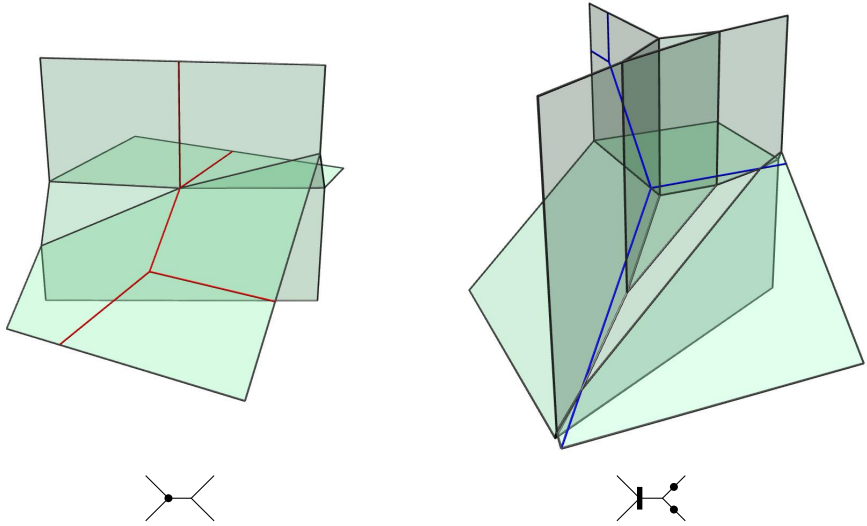

We see that the first two families have the same combinatorial type: One of the vertices lies on a vertex of the surface, while the other one is allowed to move on a halfline; see Figure \ref{fig_families}(left). These lines can in fact be explained away using an argument from \cite{BrugalleShaw}: They are not relatively realizable.

The third and last family, however, is somewhat baffling: One of the vertices lies on an \emph{edge} of the surface, while the other one can move (see Figure \ref{fig_families}, right hand side). In fact, this combinatorial type is missing from the table provided by Vigeland in \cite[Table 2]{Vigeland}. There is currently no known obstruction (i.e., non-realizability result) for any of the lines in this family. In fact, we will shortly see that even in terms of intersection combinatorics, any one of them --- even the degenerate one --- would fit the bill. Note that the vertex on the edge has coordinates $(0,0,0,0)$, so if we shrink the bounded edge to length 0, the line is actually the Bergman fan of $U_{2,4}$; see also Listing \ref{lstdegenerate}. 

\subsection{Intersection products in smooth surfaces}

We want to apply the algorithm we saw in Listing \ref{lstdegenerate} to the case of lines in a cubic surface and check that the 26 lines, together with a representative of the \enquote{odd} family we found above, fulfill the following criteria (we denote by $\cdot_V$ the intersection product in the surface $V$):
\begin{enumerate}
 \item For any line $L$, there are exactly 10 other lines $L'$, such that $L \cdot_V L' = 1$. Also, $L \cdot_V L'' = 0$ for all other lines $L''$.
 \item $L \cdot_V L = -1$ for all lines $L$.
\end{enumerate}
In fact, to save space we will only verify this here for the representative of the last family and leave it to the interested reader to complete the computation for all lines. 

\begin{lstlisting}[float,floatplacement=H,language=pmshell,caption={Intersecting lines in a smooth surface.},label=lstlineintersection]
&\shelltrop& $family_line = $rp[2];
&\shelltrop& @products = map { 
  intersect_in_smooth_surface($V,$family_line, $_)->DEGREE } 
  $L->all_isolated();
&\shelltrop& print join(",",@products);
0,1,0,0,0,0,0,0,0,0,0,1,1,1,0,0,0,0,1,1,1,1,0,1,1,0
&\shelltrop& print intersect_in_smooth_surface(
  $V,$family_line,$family_line)->DEGREE;
-1
\end{lstlisting}

One can also check that any two disjoint lines (by which we mean lines with intersection product 0) intersect five other lines and that changing the representative of the family does not affect any of the intersection multiplicities. In fact, we already saw in Listing \ref{lstdegenerate} that the degenerate representative of the last family has self-intersection -1.

\subsection{Rings of matroids}

Arbitrary intersections of tropical cycles are generally costly to compute, since they involve numerous convex hull computations.
It is therefore desirable to make use of additional information whenever possible.
One such case is the stable intersection of two tropical linear spaces.
Speyer described this in purely combinatorial terms using the underlying valuated matroids \cite{speyer}.
In the case of trivial valuation, i.e.\ when intersecting two Bergman fans, this corresponds to the operation of \emph{matroid intersection}.
In \cite{hampematroidring}, the tropical cycle ring $\mathbb{M}_n$ generated by Bergman fans of loopfree matroids on $n$ elements is studied in detail. It is shown that \emph{nested matroids} form a basis for this space and that it is, in fact, the cohomology ring of a toric variety.

\begin{theorem}
 The ring $\mathbb{M}_n$ is isomorphic to the cohomology ring $A^*(X(\textnormal{Perm}_n))$ of the toric variety corresponding to the normal fan of the permutahedron of order~$n$.
\end{theorem}

Using the explicit representation in terms of nested matroids, one can much more easily compute sums and products of cycles lying in this ring. Various properties can also be read off of this data --- such as the degree, which is simply the sum of the coefficients. It was also shown in \cite{hampematroidring}, that various matroid invariants, such as the Tutte polynomial, are linear maps on $\mathbb{M}_n$.

In the example in Listing \ref{lstmatroidring}, we consider the direct sum of two uniform matroids $U_{1,2}$. In $\mathbb{M}_4$, it is the sum of three nested matroids and we compute both its Tutte polynomial and the corresponding linear combination of the Tutte polynomials of the nested matroids to see that they are equal.

\begin{lstlisting}[language=pmshell,caption={Computing the Tutte polynomial of a direct sum of matroids.},label=lstmatroidring]
&\shelltrop& $u = matroid::uniform_matroid(1,2);
&\shelltrop& $m = matroid::direct_sum($u,$u);
&\shelltrop& print $m->TUTTE_POLYNOMIAL;
x^2 + 2*x*y + y^2
&\shelltrop& $r = matroid_ring_cycle<Min>($m);
&\shelltrop& print $r->NESTED_COEFFICIENTS;
-1 1 1
&\shelltrop& @n = $r->nested_matroids();
&\shelltrop& print - $n[0]->TUTTE_POLYNOMIAL
  + $n[1]->TUTTE_POLYNOMIAL + $n[2]->TUTTE_POLYNOMIAL ;
x^2 + 2*x*y + y^2
\end{lstlisting}

\begin{acknowledgement}
 We would like to thank Diane Maclagan for many helpful suggestions on improving this paper.
\end{acknowledgement}

\bibliographystyle{amsplain}
\bibliography{main}

\end{document}

%% file: smoothKleinTriangulation.tikz

\begin{tikzpicture}[x  = {(-0.866cm,-0.5cm)},
                    y  = {(0.866cm,-0.5cm)},
                    z  = {(0,1cm)},
                    scale = 0.75,
                    color = {lightgray}]

  \definecolor{pointcolor}{rgb}{ 1,0,0 }
  \tikzstyle{pointstyle} = [fill=pointcolor]

  \coordinate (v0) at (0, 0, 4);
  \coordinate (v1) at (1, 0, 3);
  \coordinate (v2) at (0, 1, 3);
  \coordinate (v3) at (2, 0, 2);
  \coordinate (v4) at (1, 1, 2);
  \coordinate (v5) at (0, 2, 2);
  \coordinate (v6) at (3, 0, 1);
  \coordinate (v7) at (2, 1, 1);
  \coordinate (v8) at (1, 2, 1);
  \coordinate (v9) at (0, 3, 1);
  \coordinate (v10) at (4, 0, 0);
  \coordinate (v11) at (3, 1, 0);
  \coordinate (v12) at (2, 2, 0);
  \coordinate (v13) at (1, 3, 0);
  \coordinate (v14) at (0, 4, 0);

  \definecolor{linecolor}{rgb}{ 0 0 0 }
  \tikzstyle{linestyle} = [color=linecolor, thick]

  \draw[linestyle] (v1) -- (v0);
  \draw[linestyle] (v2) -- (v0);

  \draw[linestyle] (v2) -- (v1);
  \draw[linestyle] (v3) -- (v1);
  \draw[linestyle] (v4) -- (v1);
  \draw[linestyle] (v5) -- (v1);
  \draw[linestyle] (v5) -- (v2);

  \draw[linestyle] (v6) -- (v3);
  \draw[linestyle] (v7) -- (v1);
  \draw[linestyle] (v7) -- (v4);
  \draw[linestyle] (v8) -- (v4);
  \draw[linestyle] (v8) -- (v7);
  \draw[linestyle] (v9) -- (v1);
  \draw[linestyle] (v9) -- (v4);

  \draw[linestyle] (v9) -- (v5);

  \draw[linestyle] (v9) -- (v8);
  \draw[linestyle] (v10) -- (v6);
  \draw[linestyle] (v11) -- (v1);

  \draw[linestyle] (v11) -- (v3);

  \draw[linestyle] (v11) -- (v6);

  \draw[linestyle] (v11) -- (v7);

  \draw[linestyle] (v11) -- (v8);

  \draw[linestyle] (v11) -- (v9);
  \draw[linestyle] (v11) -- (v10);

  \draw[linestyle] (v12) -- (v9);
  \draw[linestyle] (v12) -- (v11);

  \draw[linestyle] (v13) -- (v9);
  \draw[linestyle] (v13) -- (v12);

  \draw[linestyle] (v14) -- (v9);

  \draw[linestyle] (v14) -- (v13);

  \definecolor{pointcolor}{rgb}{ 1,0,0 }
  \definecolor{pointcolor}{rgb}{ 1,0,0 }
  \definecolor{pointcolor}{rgb}{ 1,0,0 }

  \coordinate (v0_unnamed__1) at (0, 0, 4);
  \coordinate (v1_unnamed__1) at (1, 0, 3);
  \coordinate (v2_unnamed__1) at (0, 1, 3);

  \definecolor{edgecolor_unnamed__1}{rgb}{ 0,0,0 }

  \definecolor{facetcolor_unnamed__1}{rgb}{ 0.4667,0.9255,0.6196 }

  \tikzstyle{facestyle_unnamed__1} = [fill=facetcolor_unnamed__1, fill opacity=0.5, draw=edgecolor_unnamed__1, line width=1 pt, line cap=round, line join=round]

  \draw[facestyle_unnamed__1] (v0_unnamed__1) -- (v1_unnamed__1) -- (v2_unnamed__1) -- (v0_unnamed__1) -- cycle;

  \fill[pointcolor] (v0_unnamed__1) circle (1 pt);
  \fill[pointcolor] (v1_unnamed__1) circle (1 pt);
  \fill[pointcolor] (v2_unnamed__1) circle (1 pt);


  \definecolor{pointcolor}{rgb}{ 1,0,0 }
  \definecolor{pointcolor}{rgb}{ 1,0,0 }
  \definecolor{pointcolor}{rgb}{ 1,0,0 }

  \coordinate (v0_unnamed__2) at (0, 3, 1);
  \coordinate (v1_unnamed__2) at (3, 1, 0);
  \coordinate (v2_unnamed__2) at (2, 2, 0);

  \definecolor{edgecolor_unnamed__2}{rgb}{ 0,0,0 }

  \definecolor{facetcolor_unnamed__2}{rgb}{ 0.4667,0.9255,0.6196 }

  \tikzstyle{facestyle_unnamed__2} = [fill=facetcolor_unnamed__2, fill opacity=0.5, draw=edgecolor_unnamed__2, line width=1 pt, line cap=round, line join=round]

  \draw[facestyle_unnamed__2] (v0_unnamed__2) -- (v1_unnamed__2) -- (v2_unnamed__2) -- (v0_unnamed__2) -- cycle;

  \fill[pointcolor] (v0_unnamed__2) circle (1 pt);
  \fill[pointcolor] (v1_unnamed__2) circle (1 pt);
  \fill[pointcolor] (v2_unnamed__2) circle (1 pt);


  \definecolor{pointcolor}{rgb}{ 1,0,0 }
  \definecolor{pointcolor}{rgb}{ 1,0,0 }
  \definecolor{pointcolor}{rgb}{ 1,0,0 }

  \coordinate (v0_unnamed__3) at (0, 3, 1);
  \coordinate (v1_unnamed__3) at (2, 2, 0);
  \coordinate (v2_unnamed__3) at (1, 3, 0);

  \definecolor{edgecolor_unnamed__3}{rgb}{ 0,0,0 }

  \definecolor{facetcolor_unnamed__3}{rgb}{ 0.4667,0.9255,0.6196 }

  \tikzstyle{facestyle_unnamed__3} = [fill=facetcolor_unnamed__3, fill opacity=0.5, draw=edgecolor_unnamed__3, line width=1 pt, line cap=round, line join=round]

  \draw[facestyle_unnamed__3] (v0_unnamed__3) -- (v1_unnamed__3) -- (v2_unnamed__3) -- (v0_unnamed__3) -- cycle;

  \fill[pointcolor] (v0_unnamed__3) circle (1 pt);
  \fill[pointcolor] (v1_unnamed__3) circle (1 pt);
  \fill[pointcolor] (v2_unnamed__3) circle (1 pt);


  \definecolor{pointcolor}{rgb}{ 1,0,0 }
  \definecolor{pointcolor}{rgb}{ 1,0,0 }
  \definecolor{pointcolor}{rgb}{ 1,0,0 }

  \coordinate (v0_unnamed__4) at (0, 3, 1);
  \coordinate (v1_unnamed__4) at (1, 3, 0);
  \coordinate (v2_unnamed__4) at (0, 4, 0);

  \definecolor{edgecolor_unnamed__4}{rgb}{ 0,0,0 }

  \definecolor{facetcolor_unnamed__4}{rgb}{ 0.4667,0.9255,0.6196 }

  \tikzstyle{facestyle_unnamed__4} = [fill=facetcolor_unnamed__4, fill opacity=0.5, draw=edgecolor_unnamed__4, line width=1 pt, line cap=round, line join=round]

  \draw[facestyle_unnamed__4] (v0_unnamed__4) -- (v1_unnamed__4) -- (v2_unnamed__4) -- (v0_unnamed__4) -- cycle;

  \fill[pointcolor] (v0_unnamed__4) circle (1 pt);
  \fill[pointcolor] (v1_unnamed__4) circle (1 pt);
  \fill[pointcolor] (v2_unnamed__4) circle (1 pt);


  \definecolor{pointcolor}{rgb}{ 1,0,0 }
  \definecolor{pointcolor}{rgb}{ 1,0,0 }
  \definecolor{pointcolor}{rgb}{ 1,0,0 }

  \coordinate (v0_unnamed__5) at (1, 0, 3);
  \coordinate (v1_unnamed__5) at (0, 1, 3);
  \coordinate (v2_unnamed__5) at (0, 2, 2);

  \definecolor{edgecolor_unnamed__5}{rgb}{ 0,0,0 }

  \definecolor{facetcolor_unnamed__5}{rgb}{ 0.4667,0.9255,0.6196 }

  \tikzstyle{facestyle_unnamed__5} = [fill=facetcolor_unnamed__5, fill opacity=0.5, draw=edgecolor_unnamed__5, line width=1 pt, line cap=round, line join=round]

  \draw[facestyle_unnamed__5] (v1_unnamed__5) -- (v0_unnamed__5) -- (v2_unnamed__5) -- (v1_unnamed__5) -- cycle;

  \fill[pointcolor] (v1_unnamed__5) circle (1 pt);
  \fill[pointcolor] (v0_unnamed__5) circle (1 pt);
  \fill[pointcolor] (v2_unnamed__5) circle (1 pt);


  \definecolor{pointcolor}{rgb}{ 1,0,0 }
  \definecolor{pointcolor}{rgb}{ 1,0,0 }
  \definecolor{pointcolor}{rgb}{ 1,0,0 }

  \coordinate (v0_unnamed__6) at (3, 0, 1);
  \coordinate (v1_unnamed__6) at (4, 0, 0);
  \coordinate (v2_unnamed__6) at (3, 1, 0);

  \definecolor{edgecolor_unnamed__6}{rgb}{ 0,0,0 }

  \definecolor{facetcolor_unnamed__6}{rgb}{ 0.4667,0.9255,0.6196 }

  \tikzstyle{facestyle_unnamed__6} = [fill=facetcolor_unnamed__6, fill opacity=0.5, draw=edgecolor_unnamed__6, line width=1 pt, line cap=round, line join=round]

  \draw[facestyle_unnamed__6] (v0_unnamed__6) -- (v1_unnamed__6) -- (v2_unnamed__6) -- (v0_unnamed__6) -- cycle;

  \fill[pointcolor] (v0_unnamed__6) circle (1 pt);
  \fill[pointcolor] (v1_unnamed__6) circle (1 pt);
  \fill[pointcolor] (v2_unnamed__6) circle (1 pt);


  \definecolor{pointcolor}{rgb}{ 1,0,0 }
  \definecolor{pointcolor}{rgb}{ 1,0,0 }
  \definecolor{pointcolor}{rgb}{ 1,0,0 }

  \coordinate (v0_unnamed__7) at (2, 0, 2);
  \coordinate (v1_unnamed__7) at (3, 0, 1);
  \coordinate (v2_unnamed__7) at (3, 1, 0);

  \definecolor{edgecolor_unnamed__7}{rgb}{ 0,0,0 }

  \definecolor{facetcolor_unnamed__7}{rgb}{ 0.4667,0.9255,0.6196 }

  \tikzstyle{facestyle_unnamed__7} = [fill=facetcolor_unnamed__7, fill opacity=0.5, draw=edgecolor_unnamed__7, line width=1 pt, line cap=round, line join=round]

  \draw[facestyle_unnamed__7] (v0_unnamed__7) -- (v2_unnamed__7) -- (v1_unnamed__7) -- (v0_unnamed__7) -- cycle;

  \fill[pointcolor] (v0_unnamed__7) circle (1 pt);
  \fill[pointcolor] (v2_unnamed__7) circle (1 pt);
  \fill[pointcolor] (v1_unnamed__7) circle (1 pt);


  \definecolor{pointcolor}{rgb}{ 1,0,0 }
  \definecolor{pointcolor}{rgb}{ 1,0,0 }
  \definecolor{pointcolor}{rgb}{ 1,0,0 }

  \coordinate (v0_unnamed__8) at (1, 0, 3);
  \coordinate (v1_unnamed__8) at (0, 2, 2);
  \coordinate (v2_unnamed__8) at (0, 3, 1);

  \definecolor{edgecolor_unnamed__8}{rgb}{ 0,0,0 }

  \definecolor{facetcolor_unnamed__8}{rgb}{ 0.4667,0.9255,0.6196 }

  \tikzstyle{facestyle_unnamed__8} = [fill=facetcolor_unnamed__8, fill opacity=0.5, draw=edgecolor_unnamed__8, line width=1 pt, line cap=round, line join=round]

  \draw[facestyle_unnamed__8] (v1_unnamed__8) -- (v0_unnamed__8) -- (v2_unnamed__8) -- (v1_unnamed__8) -- cycle;

  \fill[pointcolor] (v1_unnamed__8) circle (1 pt);
  \fill[pointcolor] (v0_unnamed__8) circle (1 pt);
  \fill[pointcolor] (v2_unnamed__8) circle (1 pt);


  \definecolor{pointcolor}{rgb}{ 1,0,0 }
  \definecolor{pointcolor}{rgb}{ 1,0,0 }
  \definecolor{pointcolor}{rgb}{ 1,0,0 }

  \coordinate (v0_unnamed__9) at (1, 0, 3);
  \coordinate (v1_unnamed__9) at (2, 0, 2);
  \coordinate (v2_unnamed__9) at (3, 1, 0);

  \definecolor{edgecolor_unnamed__9}{rgb}{ 0,0,0 }

  \definecolor{facetcolor_unnamed__9}{rgb}{ 0.4667,0.9255,0.6196 }

  \tikzstyle{facestyle_unnamed__9} = [fill=facetcolor_unnamed__9, fill opacity=0.5, draw=edgecolor_unnamed__9, line width=1 pt, line cap=round, line join=round]

  \draw[facestyle_unnamed__9] (v0_unnamed__9) -- (v2_unnamed__9) -- (v1_unnamed__9) -- (v0_unnamed__9) -- cycle;

  \fill[pointcolor] (v0_unnamed__9) circle (1 pt);
  \fill[pointcolor] (v2_unnamed__9) circle (1 pt);
  \fill[pointcolor] (v1_unnamed__9) circle (1 pt);


  \definecolor{pointcolor}{rgb}{ 1,0,0 }
  \definecolor{pointcolor}{rgb}{ 1,0,0 }
  \definecolor{pointcolor}{rgb}{ 1,0,0 }

  \coordinate (v0_unnamed__10) at (1, 2, 1);
  \coordinate (v1_unnamed__10) at (0, 3, 1);
  \coordinate (v2_unnamed__10) at (3, 1, 0);

  \definecolor{edgecolor_unnamed__10}{rgb}{ 0,0,0 }

  \definecolor{facetcolor_unnamed__10}{rgb}{ 0.4667,0.9255,0.6196 }

  \tikzstyle{facestyle_unnamed__10} = [fill=facetcolor_unnamed__10, fill opacity=0.5, draw=edgecolor_unnamed__10, line width=1 pt, line cap=round, line join=round]

  \draw[facestyle_unnamed__10] (v0_unnamed__10) -- (v2_unnamed__10) -- (v1_unnamed__10) -- (v0_unnamed__10) -- cycle;

  \fill[pointcolor] (v0_unnamed__10) circle (1 pt);
  \fill[pointcolor] (v2_unnamed__10) circle (1 pt);
  \fill[pointcolor] (v1_unnamed__10) circle (1 pt);


  \definecolor{pointcolor}{rgb}{ 1,0,0 }
  \definecolor{pointcolor}{rgb}{ 1,0,0 }
  \definecolor{pointcolor}{rgb}{ 1,0,0 }

  \coordinate (v0_unnamed__11) at (1, 0, 3);
  \coordinate (v1_unnamed__11) at (1, 1, 2);
  \coordinate (v2_unnamed__11) at (0, 3, 1);

  \definecolor{edgecolor_unnamed__11}{rgb}{ 0,0,0 }

  \definecolor{facetcolor_unnamed__11}{rgb}{ 0.4667,0.9255,0.6196 }

  \tikzstyle{facestyle_unnamed__11} = [fill=facetcolor_unnamed__11, fill opacity=0.5, draw=edgecolor_unnamed__11, line width=1 pt, line cap=round, line join=round]

  \draw[facestyle_unnamed__11] (v0_unnamed__11) -- (v1_unnamed__11) -- (v2_unnamed__11) -- (v0_unnamed__11) -- cycle;

  \fill[pointcolor] (v0_unnamed__11) circle (1 pt);
  \fill[pointcolor] (v1_unnamed__11) circle (1 pt);
  \fill[pointcolor] (v2_unnamed__11) circle (1 pt);


  \definecolor{pointcolor}{rgb}{ 1,0,0 }
  \definecolor{pointcolor}{rgb}{ 1,0,0 }
  \definecolor{pointcolor}{rgb}{ 1,0,0 }

  \coordinate (v0_unnamed__12) at (1, 0, 3);
  \coordinate (v1_unnamed__12) at (2, 1, 1);
  \coordinate (v2_unnamed__12) at (3, 1, 0);

  \definecolor{edgecolor_unnamed__12}{rgb}{ 0,0,0 }

  \definecolor{facetcolor_unnamed__12}{rgb}{ 0.4667,0.9255,0.6196 }

  \tikzstyle{facestyle_unnamed__12} = [fill=facetcolor_unnamed__12, fill opacity=0.5, draw=edgecolor_unnamed__12, line width=1 pt, line cap=round, line join=round]

  \draw[facestyle_unnamed__12] (v0_unnamed__12) -- (v2_unnamed__12) -- (v1_unnamed__12) -- (v0_unnamed__12) -- cycle;

  \fill[pointcolor] (v0_unnamed__12) circle (1 pt);
  \fill[pointcolor] (v2_unnamed__12) circle (1 pt);
  \fill[pointcolor] (v1_unnamed__12) circle (1 pt);


  \definecolor{pointcolor}{rgb}{ 1,0,0 }
  \definecolor{pointcolor}{rgb}{ 1,0,0 }
  \definecolor{pointcolor}{rgb}{ 1,0,0 }

  \coordinate (v0_unnamed__13) at (2, 1, 1);
  \coordinate (v1_unnamed__13) at (1, 2, 1);
  \coordinate (v2_unnamed__13) at (3, 1, 0);

  \definecolor{edgecolor_unnamed__13}{rgb}{ 0,0,0 }

  \definecolor{facetcolor_unnamed__13}{rgb}{ 0.4667,0.9255,0.6196 }

  \tikzstyle{facestyle_unnamed__13} = [fill=facetcolor_unnamed__13, fill opacity=0.5, draw=edgecolor_unnamed__13, line width=1 pt, line cap=round, line join=round]

  \draw[facestyle_unnamed__13] (v0_unnamed__13) -- (v2_unnamed__13) -- (v1_unnamed__13) -- (v0_unnamed__13) -- cycle;

  \fill[pointcolor] (v0_unnamed__13) circle (1 pt);
  \fill[pointcolor] (v2_unnamed__13) circle (1 pt);
  \fill[pointcolor] (v1_unnamed__13) circle (1 pt);


  \definecolor{pointcolor}{rgb}{ 1,0,0 }
  \definecolor{pointcolor}{rgb}{ 1,0,0 }
  \definecolor{pointcolor}{rgb}{ 1,0,0 }

  \coordinate (v0_unnamed__14) at (1, 1, 2);
  \coordinate (v1_unnamed__14) at (1, 2, 1);
  \coordinate (v2_unnamed__14) at (0, 3, 1);

  \definecolor{edgecolor_unnamed__14}{rgb}{ 0,0,0 }

  \definecolor{facetcolor_unnamed__14}{rgb}{ 0.4667,0.9255,0.6196 }

  \tikzstyle{facestyle_unnamed__14} = [fill=facetcolor_unnamed__14, fill opacity=0.5, draw=edgecolor_unnamed__14, line width=1 pt, line cap=round, line join=round]

  \draw[facestyle_unnamed__14] (v0_unnamed__14) -- (v1_unnamed__14) -- (v2_unnamed__14) -- (v0_unnamed__14) -- cycle;

  \fill[pointcolor] (v0_unnamed__14) circle (1 pt);
  \fill[pointcolor] (v1_unnamed__14) circle (1 pt);
  \fill[pointcolor] (v2_unnamed__14) circle (1 pt);


  \definecolor{pointcolor}{rgb}{ 1,0,0 }
  \definecolor{pointcolor}{rgb}{ 1,0,0 }
  \definecolor{pointcolor}{rgb}{ 1,0,0 }

  \coordinate (v0_unnamed__15) at (1, 0, 3);
  \coordinate (v1_unnamed__15) at (1, 1, 2);
  \coordinate (v2_unnamed__15) at (2, 1, 1);

  \definecolor{edgecolor_unnamed__15}{rgb}{ 0,0,0 }

  \definecolor{facetcolor_unnamed__15}{rgb}{ 0.4667,0.9255,0.6196 }

  \tikzstyle{facestyle_unnamed__15} = [fill=facetcolor_unnamed__15, fill opacity=0.5, draw=edgecolor_unnamed__15, line width=1 pt, line cap=round, line join=round]

  \draw[facestyle_unnamed__15] (v0_unnamed__15) -- (v2_unnamed__15) -- (v1_unnamed__15) -- (v0_unnamed__15) -- cycle;

  \fill[pointcolor] (v0_unnamed__15) circle (1 pt);
  \fill[pointcolor] (v2_unnamed__15) circle (1 pt);
  \fill[pointcolor] (v1_unnamed__15) circle (1 pt);


  \definecolor{pointcolor}{rgb}{ 1,0,0 }
  \definecolor{pointcolor}{rgb}{ 1,0,0 }
  \definecolor{pointcolor}{rgb}{ 1,0,0 }

  \coordinate (v0_unnamed__16) at (1, 1, 2);
  \coordinate (v1_unnamed__16) at (2, 1, 1);
  \coordinate (v2_unnamed__16) at (1, 2, 1);

  \definecolor{edgecolor_unnamed__16}{rgb}{ 0,0,0 }

  \definecolor{facetcolor_unnamed__16}{rgb}{ 0.4667,0.9255,0.6196 }

  \tikzstyle{facestyle_unnamed__16} = [fill=facetcolor_unnamed__16, fill opacity=0.5, draw=edgecolor_unnamed__16, line width=1 pt, line cap=round, line join=round]

  \draw[facestyle_unnamed__16] (v0_unnamed__16) -- (v1_unnamed__16) -- (v2_unnamed__16) -- (v0_unnamed__16) -- cycle;

  \fill[pointcolor] (v0_unnamed__16) circle (1 pt);
  \fill[pointcolor] (v1_unnamed__16) circle (1 pt);
  \fill[pointcolor] (v2_unnamed__16) circle (1 pt);


  \fill[pointcolor] (v0) circle (1 pt);
  \draw (v0) node [draw,circle,text=black, inner sep=0pt, minimum size=15pt, fill=white,align=left] {0};

  \fill[pointcolor] (v1) circle (1 pt);
  \node at (v1) [draw,circle,text=black, inner sep=0pt, minimum size=15pt, fill=white,align=left] {1};

  \fill[pointcolor] (v2) circle (1 pt);
  \node at (v2) [draw,circle,text=black, inner sep=0pt, minimum size=15pt, fill=white,align=left] {2};

  \fill[pointcolor] (v3) circle (1 pt);
  \node at (v3) [draw,circle,text=black, inner sep=0pt, minimum size=15pt, fill=white,align=left] {3};

  \fill[pointcolor] (v4) circle (1 pt);
  \node at (v4) [draw,circle,text=black, inner sep=0pt, minimum size=15pt, fill=white,align=left] {4};

  \fill[pointcolor] (v5) circle (1 pt);
  \node at (v5) [draw,circle,text=black, inner sep=0pt, minimum size=15pt, fill=white,align=left] {5};

  \fill[pointcolor] (v6) circle (1 pt);
  \node at (v6) [draw,circle,text=black, inner sep=0pt, minimum size=15pt, fill=white,align=left] {6};

  \fill[pointcolor] (v7) circle (1 pt);
  \node at (v7) [draw,circle,text=black, inner sep=0pt, minimum size=15pt, fill=white,align=left] {7};

  \fill[pointcolor] (v8) circle (1 pt);
  \node at (v8) [draw,circle,text=black, inner sep=0pt, minimum size=15pt, fill=white,align=left] {8};

  \fill[pointcolor] (v9) circle (1 pt);
  \node at (v9) [draw,circle,text=black, inner sep=0pt, minimum size=15pt, fill=white,align=left] {9};

  \fill[pointcolor] (v10) circle (1 pt);
  \node at (v10) [draw,circle,text=black, inner sep=0pt, minimum size=15pt, fill=white,align=left] {10};

  \fill[pointcolor] (v11) circle (1 pt);
  \node at (v11) [draw,circle,text=black, inner sep=0pt, minimum size=15pt, fill=white,align=left] {11};

  \fill[pointcolor] (v12) circle (1 pt);
  \node at (v12) [draw,circle,text=black, inner sep=0pt, minimum size=15pt, fill=white,align=left] {12};

  \fill[pointcolor] (v13) circle (1 pt);
  \node at (v13) [draw,circle,text=black, inner sep=0pt, minimum size=15pt, fill=white,align=left] {13};

  \fill[pointcolor] (v14) circle (1 pt);
  \node at (v14) [draw,circle,text=black, inner sep=0pt, minimum size=15pt, fill=white,align=left] {14};

\end{tikzpicture}


%% file: smoothKleinTriangulation-ray0.tikz

\begin{tikzpicture}[x  = {(-0.866cm,-0.5cm)},
                    y  = {(0.866cm,-0.5cm)},
                    z  = {(0,1cm)},
                    scale = 0.4,
                    color = {lightgray}]

  \definecolor{pointcolor}{rgb}{ 0,0,0 }
  \tikzstyle{pointstyle} = [fill=pointcolor]

  \coordinate (v0) at (0, 0, 4);
  \coordinate (v1) at (1, 0, 3);
  \coordinate (v2) at (0, 1, 3);
  \coordinate (v3) at (2, 0, 2);
  \coordinate (v4) at (1, 1, 2);
  \coordinate (v5) at (0, 2, 2);
  \coordinate (v6) at (3, 0, 1);
  \coordinate (v7) at (2, 1, 1);
  \coordinate (v8) at (1, 2, 1);
  \coordinate (v9) at (0, 3, 1);
  \coordinate (v10) at (4, 0, 0);
  \coordinate (v11) at (3, 1, 0);
  \coordinate (v12) at (2, 2, 0);
  \coordinate (v13) at (1, 3, 0);
  \coordinate (v14) at (0, 4, 0);

  \definecolor{linecolor}{rgb}{ 0 0 0 }
  \tikzstyle{linestyle} = [color=linecolor, thick]

  \definecolor{facetcolor}{rgb}{ 0.4667,0.9255,0.6196 }
  \tikzstyle{facetstyle} = [fill=facetcolor, fill opacity=0.5, draw=linecolor, line width=1 pt, line cap=round, line join=round]

  \draw[facetstyle] (v0) -- (v10) -- (v14) -- cycle;
  \draw[linestyle] (v1) -- (v2);
  
  \fill[white,draw=black] (v0) circle (4 pt);
  \node at (v0) [text=black, inner sep=0.5pt, above right, draw=none, align=left] {};

  \fill[white,draw=black] (v1) circle (4 pt);
  \node at (v1) [text=black, inner sep=0.5pt, above right, draw=none, align=left] {};

  \fill[white,draw=black] (v2) circle (4 pt);
  \node at (v2) [text=black, inner sep=0.5pt, above right, draw=none, align=left] {};

  \fill[white,draw=black] (v3) circle (4 pt);
  \node at (v3) [text=black, inner sep=0.5pt, above right, draw=none, align=left] {};

  \fill[white,draw=black] (v4) circle (4 pt);
  \node at (v4) [text=black, inner sep=0.5pt, above right, draw=none, align=left] {};

  \fill[white,draw=black] (v5) circle (4 pt);
  \node at (v5) [text=black, inner sep=0.5pt, above right, draw=none, align=left] {};

  \fill[white,draw=black] (v6) circle (4 pt);
  \node at (v6) [text=black, inner sep=0.5pt, above right, draw=none, align=left] {};

  \fill[white,draw=black] (v7) circle (4 pt);
  \node at (v7) [text=black, inner sep=0.5pt, above right, draw=none, align=left] {};

  \fill[white,draw=black] (v8) circle (4 pt);
  \node at (v8) [text=black, inner sep=0.5pt, above right, draw=none, align=left] {};

  \fill[white,draw=black] (v9) circle (4 pt);
  \node at (v9) [text=black, inner sep=0.5pt, above right, draw=none, align=left] {};

  \fill[white,draw=black] (v10) circle (4 pt);
  \node at (v10) [text=black, inner sep=0.5pt, above right, draw=none, align=left] {};

  \fill[white,draw=black] (v11) circle (4 pt);
  \node at (v11) [text=black, inner sep=0.5pt, above right, draw=none, align=left] {};

  \fill[white,draw=black] (v12) circle (4 pt);
  \node at (v12) [text=black, inner sep=0.5pt, above right, draw=none, align=left] {};

  \fill[white,draw=black] (v13) circle (4 pt);
  \node at (v13) [text=black, inner sep=0.5pt, above right, draw=none, align=left] {};

  \fill[white,draw=black] (v14) circle (4 pt);
  \node at (v14) [text=black, inner sep=0.5pt, above right, draw=none, align=left] {};

\end{tikzpicture}

%% file: smoothKleinTriangulation-ray2.tikz

\begin{tikzpicture}[x  = {(-0.866cm,-0.5cm)},
                    y  = {(0.866cm,-0.5cm)},
                    z  = {(0,1cm)},
                    scale = 0.4,
                    color = {lightgray}]

  \definecolor{pointcolor}{rgb}{ 0,0,0 }
  \tikzstyle{pointstyle} = [fill=pointcolor]

  \coordinate (v0) at (0, 0, 4);
  \coordinate (v1) at (1, 0, 3);
  \coordinate (v2) at (0, 1, 3);
  \coordinate (v3) at (2, 0, 2);
  \coordinate (v4) at (1, 1, 2);
  \coordinate (v5) at (0, 2, 2);
  \coordinate (v6) at (3, 0, 1);
  \coordinate (v7) at (2, 1, 1);
  \coordinate (v8) at (1, 2, 1);
  \coordinate (v9) at (0, 3, 1);
  \coordinate (v10) at (4, 0, 0);
  \coordinate (v11) at (3, 1, 0);
  \coordinate (v12) at (2, 2, 0);
  \coordinate (v13) at (1, 3, 0);
  \coordinate (v14) at (0, 4, 0);

  \definecolor{linecolor}{rgb}{ 0 0 0 }
  \tikzstyle{linestyle} = [color=linecolor, thick]

  \definecolor{facetcolor}{rgb}{ 0.4667,0.9255,0.6196 }
  \tikzstyle{facetstyle} = [fill=facetcolor, fill opacity=0.5, draw=linecolor, line width=1 pt, line cap=round, line join=round]

  \draw[facetstyle] (v0) -- (v10) -- (v14) -- cycle;
  \draw[linestyle] (v1) -- (v5);
  
  \fill[white,draw=black] (v0) circle (4 pt);
  \node at (v0) [text=black, inner sep=0.5pt, above right, draw=none, align=left] {};

  \fill[white,draw=black] (v1) circle (4 pt);
  \node at (v1) [text=black, inner sep=0.5pt, above right, draw=none, align=left] {};

  \fill[white,draw=black] (v2) circle (4 pt);
  \node at (v2) [text=black, inner sep=0.5pt, above right, draw=none, align=left] {};

  \fill[white,draw=black] (v3) circle (4 pt);
  \node at (v3) [text=black, inner sep=0.5pt, above right, draw=none, align=left] {};

  \fill[white,draw=black] (v4) circle (4 pt);
  \node at (v4) [text=black, inner sep=0.5pt, above right, draw=none, align=left] {};

  \fill[white,draw=black] (v5) circle (4 pt);
  \node at (v5) [text=black, inner sep=0.5pt, above right, draw=none, align=left] {};

  \fill[white,draw=black] (v6) circle (4 pt);
  \node at (v6) [text=black, inner sep=0.5pt, above right, draw=none, align=left] {};

  \fill[white,draw=black] (v7) circle (4 pt);
  \node at (v7) [text=black, inner sep=0.5pt, above right, draw=none, align=left] {};

  \fill[white,draw=black] (v8) circle (4 pt);
  \node at (v8) [text=black, inner sep=0.5pt, above right, draw=none, align=left] {};

  \fill[white,draw=black] (v9) circle (4 pt);
  \node at (v9) [text=black, inner sep=0.5pt, above right, draw=none, align=left] {};

  \fill[white,draw=black] (v10) circle (4 pt);
  \node at (v10) [text=black, inner sep=0.5pt, above right, draw=none, align=left] {};

  \fill[white,draw=black] (v11) circle (4 pt);
  \node at (v11) [text=black, inner sep=0.5pt, above right, draw=none, align=left] {};

  \fill[white,draw=black] (v12) circle (4 pt);
  \node at (v12) [text=black, inner sep=0.5pt, above right, draw=none, align=left] {};

  \fill[white,draw=black] (v13) circle (4 pt);
  \node at (v13) [text=black, inner sep=0.5pt, above right, draw=none, align=left] {};

  \fill[white,draw=black] (v14) circle (4 pt);
  \node at (v14) [text=black, inner sep=0.5pt, above right, draw=none, align=left] {};

\end{tikzpicture}

%% file: smoothKleinTriangulation-ray8.tikz

\begin{tikzpicture}[x  = {(-0.866cm,-0.5cm)},
                    y  = {(0.866cm,-0.5cm)},
                    z  = {(0,1cm)},
                    scale = 0.4,
                    color = {lightgray}]

  \definecolor{pointcolor}{rgb}{ 0,0,0 }
  \tikzstyle{pointstyle} = [fill=pointcolor]

  \coordinate (v0) at (0, 0, 4);
  \coordinate (v1) at (1, 0, 3);
  \coordinate (v2) at (0, 1, 3);
  \coordinate (v3) at (2, 0, 2);
  \coordinate (v4) at (1, 1, 2);
  \coordinate (v5) at (0, 2, 2);
  \coordinate (v6) at (3, 0, 1);
  \coordinate (v7) at (2, 1, 1);
  \coordinate (v8) at (1, 2, 1);
  \coordinate (v9) at (0, 3, 1);
  \coordinate (v10) at (4, 0, 0);
  \coordinate (v11) at (3, 1, 0);
  \coordinate (v12) at (2, 2, 0);
  \coordinate (v13) at (1, 3, 0);
  \coordinate (v14) at (0, 4, 0);
  \definecolor{linecolor}{rgb}{ 0 0 0 }
  \tikzstyle{linestyle} = [color=linecolor, thick]

  \definecolor{facetcolor}{rgb}{ 0.4667,0.9255,0.6196 }
  \tikzstyle{facetstyle} = [fill=facetcolor, fill opacity=0.5, draw=linecolor, line width=1 pt, line cap=round, line join=round]

  \draw[facetstyle] (v0) -- (v10) -- (v14) -- cycle;
  \draw[linestyle] (v1) -- (v9);
  
  \fill[white,draw=black] (v0) circle (4 pt);
  \node at (v0) [text=black, inner sep=0.5pt, above right, draw=none, align=left] {};

  \fill[white,draw=black] (v1) circle (4 pt);
  \node at (v1) [text=black, inner sep=0.5pt, above right, draw=none, align=left] {};

  \fill[white,draw=black] (v2) circle (4 pt);
  \node at (v2) [text=black, inner sep=0.5pt, above right, draw=none, align=left] {};

  \fill[white,draw=black] (v3) circle (4 pt);
  \node at (v3) [text=black, inner sep=0.5pt, above right, draw=none, align=left] {};

  \fill[white,draw=black] (v4) circle (4 pt);
  \node at (v4) [text=black, inner sep=0.5pt, above right, draw=none, align=left] {};

  \fill[white,draw=black] (v5) circle (4 pt);
  \node at (v5) [text=black, inner sep=0.5pt, above right, draw=none, align=left] {};

  \fill[white,draw=black] (v6) circle (4 pt);
  \node at (v6) [text=black, inner sep=0.5pt, above right, draw=none, align=left] {};

  \fill[white,draw=black] (v7) circle (4 pt);
  \node at (v7) [text=black, inner sep=0.5pt, above right, draw=none, align=left] {};

  \fill[white,draw=black] (v8) circle (4 pt);
  \node at (v8) [text=black, inner sep=0.5pt, above right, draw=none, align=left] {};

  \fill[white,draw=black] (v9) circle (4 pt);
  \node at (v9) [text=black, inner sep=0.5pt, above right, draw=none, align=left] {};

  \fill[white,draw=black] (v10) circle (4 pt);
  \node at (v10) [text=black, inner sep=0.5pt, above right, draw=none, align=left] {};

  \fill[white,draw=black] (v11) circle (4 pt);
  \node at (v11) [text=black, inner sep=0.5pt, above right, draw=none, align=left] {};

  \fill[white,draw=black] (v12) circle (4 pt);
  \node at (v12) [text=black, inner sep=0.5pt, above right, draw=none, align=left] {};

  \fill[white,draw=black] (v13) circle (4 pt);
  \node at (v13) [text=black, inner sep=0.5pt, above right, draw=none, align=left] {};

  \fill[white,draw=black] (v14) circle (4 pt);
  \node at (v14) [text=black, inner sep=0.5pt, above right, draw=none, align=left] {};

\end{tikzpicture}

%% file: smoothKleinTriangulation-ray11.tikz

\begin{tikzpicture}[x  = {(-0.866cm,-0.5cm)},
                    y  = {(0.866cm,-0.5cm)},
                    z  = {(0,1cm)},
                    scale = 0.4,
                    color = {lightgray}]

  \definecolor{pointcolor}{rgb}{ 0,0,0 }
  \tikzstyle{pointstyle} = [fill=pointcolor]

  \coordinate (v0) at (0, 0, 4);
  \coordinate (v1) at (1, 0, 3);
  \coordinate (v2) at (0, 1, 3);
  \coordinate (v3) at (2, 0, 2);
  \coordinate (v4) at (1, 1, 2);
  \coordinate (v5) at (0, 2, 2);
  \coordinate (v6) at (3, 0, 1);
  \coordinate (v7) at (2, 1, 1);
  \coordinate (v8) at (1, 2, 1);
  \coordinate (v9) at (0, 3, 1);
  \coordinate (v10) at (4, 0, 0);
  \coordinate (v11) at (3, 1, 0);
  \coordinate (v12) at (2, 2, 0);
  \coordinate (v13) at (1, 3, 0);
  \coordinate (v14) at (0, 4, 0);

  \definecolor{linecolor}{rgb}{ 0 0 0 }
  \tikzstyle{linestyle} = [color=linecolor, thick]

  \definecolor{facetcolor}{rgb}{ 0.4667,0.9255,0.6196 }
  \tikzstyle{facetstyle} = [fill=facetcolor, fill opacity=0.5, draw=linecolor, line width=1 pt, line cap=round, line join=round]

  \draw[facetstyle] (v0) -- (v10) -- (v14) -- cycle;
  \draw[linestyle] (v4) -- (v1);
  \draw[linestyle] (v4) -- (v9);
  \draw[linestyle] (v4) -- (v11);
  
  \fill[white,draw=black] (v0) circle (4 pt);
  \node at (v0) [text=black, inner sep=0.5pt, above right, draw=none, align=left] {};

  \fill[white,draw=black] (v1) circle (4 pt);
  \node at (v1) [text=black, inner sep=0.5pt, above right, draw=none, align=left] {};

  \fill[white,draw=black] (v2) circle (4 pt);
  \node at (v2) [text=black, inner sep=0.5pt, above right, draw=none, align=left] {};

  \fill[white,draw=black] (v3) circle (4 pt);
  \node at (v3) [text=black, inner sep=0.5pt, above right, draw=none, align=left] {};

  \fill[white,draw=black] (v4) circle (4 pt);
  \node at (v4) [text=black, inner sep=0.5pt, above right, draw=none, align=left] {};

  \fill[white,draw=black] (v5) circle (4 pt);
  \node at (v5) [text=black, inner sep=0.5pt, above right, draw=none, align=left] {};

  \fill[white,draw=black] (v6) circle (4 pt);
  \node at (v6) [text=black, inner sep=0.5pt, above right, draw=none, align=left] {};

  \fill[white,draw=black] (v7) circle (4 pt);
  \node at (v7) [text=black, inner sep=0.5pt, above right, draw=none, align=left] {};

  \fill[white,draw=black] (v8) circle (4 pt);
  \node at (v8) [text=black, inner sep=0.5pt, above right, draw=none, align=left] {};

  \fill[white,draw=black] (v9) circle (4 pt);
  \node at (v9) [text=black, inner sep=0.5pt, above right, draw=none, align=left] {};

  \fill[white,draw=black] (v10) circle (4 pt);
  \node at (v10) [text=black, inner sep=0.5pt, above right, draw=none, align=left] {};

  \fill[white,draw=black] (v11) circle (4 pt);
  \node at (v11) [text=black, inner sep=0.5pt, above right, draw=none, align=left] {};

  \fill[white,draw=black] (v12) circle (4 pt);
  \node at (v12) [text=black, inner sep=0.5pt, above right, draw=none, align=left] {};

  \fill[white,draw=black] (v13) circle (4 pt);
  \node at (v13) [text=black, inner sep=0.5pt, above right, draw=none, align=left] {};

  \fill[white,draw=black] (v14) circle (4 pt);
  \node at (v14) [text=black, inner sep=0.5pt, above right, draw=none, align=left] {};

\end{tikzpicture}

%% file: smoothKleinCurve.tikz

\begin{tikzpicture}[scale=0.75]

  \definecolor{pointcolor_unnamed__1}{rgb}{ 1,0,0 }
  \tikzstyle{pointstyle_unnamed__1} = [fill=pointcolor_unnamed__1]

  \coordinate (v0_unnamed__1) at (-3, -6);
  \coordinate (v1_unnamed__1) at (7, -6);

  \definecolor{linecolor_unnamed__1}{rgb}{ 0 0 0 }

  \tikzstyle{linestyle_unnamed__1} = [color=linecolor_unnamed__1, thick]

  \draw[linestyle_unnamed__1] (v1_unnamed__1) -- (v0_unnamed__1);

  \fill[pointcolor_unnamed__1] (v1_unnamed__1) circle (0.1 pt);
  \fill[pointcolor_unnamed__1] (v0_unnamed__1) circle (0.1 pt);


  \definecolor{pointcolor_unnamed__2}{rgb}{ 1,0,0 }
  \tikzstyle{pointstyle_unnamed__2} = [fill=pointcolor_unnamed__2]

  \coordinate (v0_unnamed__2) at (-3, -6);
  \coordinate (v1_unnamed__2) at (-4, -7);

  \definecolor{linecolor_unnamed__2}{rgb}{ 0 0 0 }

  \tikzstyle{linestyle_unnamed__2} = [color=linecolor_unnamed__2, thick]

  \draw[linestyle_unnamed__2] (v1_unnamed__2) -- (v0_unnamed__2);

  \fill[pointcolor_unnamed__2] (v1_unnamed__2) circle (0.1 pt);
  \fill[pointcolor_unnamed__2] (v0_unnamed__2) circle (0.1 pt);


  \definecolor{pointcolor_unnamed__3}{rgb}{ 1,0,0 }
  \tikzstyle{pointstyle_unnamed__3} = [fill=pointcolor_unnamed__3]

  \coordinate (v0_unnamed__3) at (-3, -5);
  \coordinate (v1_unnamed__3) at (-3, -6);

  \definecolor{linecolor_unnamed__3}{rgb}{ 0 0 0 }

  \tikzstyle{linestyle_unnamed__3} = [color=linecolor_unnamed__3, thick]

  \draw[linestyle_unnamed__3] (v1_unnamed__3) -- (v0_unnamed__3);

  \fill[pointcolor_unnamed__3] (v1_unnamed__3) circle (0.1 pt);
  \fill[pointcolor_unnamed__3] (v0_unnamed__3) circle (0.1 pt);


  \definecolor{pointcolor_unnamed__4}{rgb}{ 1,0,0 }
  \tikzstyle{pointstyle_unnamed__4} = [fill=pointcolor_unnamed__4]

  \coordinate (v0_unnamed__4) at (3, 1);
  \coordinate (v1_unnamed__4) at (7, 1);

  \definecolor{linecolor_unnamed__4}{rgb}{ 0 0 0 }

  \tikzstyle{linestyle_unnamed__4} = [color=linecolor_unnamed__4, thick]

  \draw[linestyle_unnamed__4] (v1_unnamed__4) -- (v0_unnamed__4);

  \fill[pointcolor_unnamed__4] (v1_unnamed__4) circle (0.1 pt);
  \fill[pointcolor_unnamed__4] (v0_unnamed__4) circle (0.1 pt);


  \definecolor{pointcolor_unnamed__5}{rgb}{ 1,0,0 }
  \tikzstyle{pointstyle_unnamed__5} = [fill=pointcolor_unnamed__5]

  \coordinate (v0_unnamed__5) at (0.333333, 0);
  \coordinate (v1_unnamed__5) at (-0.666667, -1);

  \definecolor{linecolor_unnamed__5}{rgb}{ 0 0 0 }

  \tikzstyle{linestyle_unnamed__5} = [color=linecolor_unnamed__5, thick]

  \draw[linestyle_unnamed__5] (v1_unnamed__5) -- (v0_unnamed__5);

  \fill[pointcolor_unnamed__5] (v1_unnamed__5) circle (0.1 pt);
  \fill[pointcolor_unnamed__5] (v0_unnamed__5) circle (0.1 pt);


  \definecolor{pointcolor_unnamed__6}{rgb}{ 1,0,0 }
  \tikzstyle{pointstyle_unnamed__6} = [fill=pointcolor_unnamed__6]

  \coordinate (v0_unnamed__6) at (-2, -3);
  \coordinate (v1_unnamed__6) at (-3, -5);

  \definecolor{linecolor_unnamed__6}{rgb}{ 0 0 0 }

  \tikzstyle{linestyle_unnamed__6} = [color=linecolor_unnamed__6, thick]

  \draw[linestyle_unnamed__6] (v1_unnamed__6) -- (v0_unnamed__6);

  \fill[pointcolor_unnamed__6] (v1_unnamed__6) circle (0.1 pt);
  \fill[pointcolor_unnamed__6] (v0_unnamed__6) circle (0.1 pt);


  \definecolor{pointcolor_unnamed__7}{rgb}{ 1,0,0 }
  \tikzstyle{pointstyle_unnamed__7} = [fill=pointcolor_unnamed__7]

  \coordinate (v0_unnamed__7) at (-3, -5);
  \coordinate (v1_unnamed__7) at (-5, -7);

  \definecolor{linecolor_unnamed__7}{rgb}{ 0 0 0 }

  \tikzstyle{linestyle_unnamed__7} = [color=linecolor_unnamed__7, thick]

  \draw[linestyle_unnamed__7] (v1_unnamed__7) -- (v0_unnamed__7);

  \fill[pointcolor_unnamed__7] (v1_unnamed__7) circle (0.1 pt);
  \fill[pointcolor_unnamed__7] (v0_unnamed__7) circle (0.1 pt);


  \definecolor{pointcolor_unnamed__8}{rgb}{ 1,0,0 }
  \tikzstyle{pointstyle_unnamed__8} = [fill=pointcolor_unnamed__8]

  \coordinate (v0_unnamed__8) at (5, 2);
  \coordinate (v1_unnamed__8) at (7, 2);

  \definecolor{linecolor_unnamed__8}{rgb}{ 0 0 0 }

  \tikzstyle{linestyle_unnamed__8} = [color=linecolor_unnamed__8, thick]

  \draw[linestyle_unnamed__8] (v1_unnamed__8) -- (v0_unnamed__8);

  \fill[pointcolor_unnamed__8] (v1_unnamed__8) circle (0.1 pt);
  \fill[pointcolor_unnamed__8] (v0_unnamed__8) circle (0.1 pt);


  \definecolor{pointcolor_unnamed__9}{rgb}{ 1,0,0 }
  \tikzstyle{pointstyle_unnamed__9} = [fill=pointcolor_unnamed__9]

  \coordinate (v0_unnamed__9) at (1, 0.333333);
  \coordinate (v1_unnamed__9) at (0.333333, 0);

  \definecolor{linecolor_unnamed__9}{rgb}{ 0 0 0 }

  \tikzstyle{linestyle_unnamed__9} = [color=linecolor_unnamed__9, thick]

  \draw[linestyle_unnamed__9] (v1_unnamed__9) -- (v0_unnamed__9);

  \fill[pointcolor_unnamed__9] (v1_unnamed__9) circle (0.1 pt);
  \fill[pointcolor_unnamed__9] (v0_unnamed__9) circle (0.1 pt);


  \definecolor{pointcolor_unnamed__10}{rgb}{ 1,0,0 }
  \tikzstyle{pointstyle_unnamed__10} = [fill=pointcolor_unnamed__10]

  \coordinate (v0_unnamed__10) at (0.333333, 0);
  \coordinate (v1_unnamed__10) at (0, 0);

  \definecolor{linecolor_unnamed__10}{rgb}{ 0 0 0 }

  \tikzstyle{linestyle_unnamed__10} = [color=linecolor_unnamed__10, thick]

  \draw[linestyle_unnamed__10] (v1_unnamed__10) -- (v0_unnamed__10);

  \fill[pointcolor_unnamed__10] (v1_unnamed__10) circle (0.1 pt);
  \fill[pointcolor_unnamed__10] (v0_unnamed__10) circle (0.1 pt);


  \definecolor{pointcolor_unnamed__11}{rgb}{ 1,0,0 }
  \tikzstyle{pointstyle_unnamed__11} = [fill=pointcolor_unnamed__11]

  \coordinate (v0_unnamed__11) at (0, 0);
  \coordinate (v1_unnamed__11) at (-0.333333, -0.333333);

  \definecolor{linecolor_unnamed__11}{rgb}{ 0 0 0 }

  \tikzstyle{linestyle_unnamed__11} = [color=linecolor_unnamed__11, thick]

  \draw[linestyle_unnamed__11] (v1_unnamed__11) -- (v0_unnamed__11);

  \fill[pointcolor_unnamed__11] (v1_unnamed__11) circle (0.1 pt);
  \fill[pointcolor_unnamed__11] (v0_unnamed__11) circle (0.1 pt);


  \definecolor{pointcolor_unnamed__12}{rgb}{ 1,0,0 }
  \tikzstyle{pointstyle_unnamed__12} = [fill=pointcolor_unnamed__12]

  \coordinate (v0_unnamed__12) at (0, 0.333333);
  \coordinate (v1_unnamed__12) at (0, 0);

  \definecolor{linecolor_unnamed__12}{rgb}{ 0 0 0 }

  \tikzstyle{linestyle_unnamed__12} = [color=linecolor_unnamed__12, thick]

  \draw[linestyle_unnamed__12] (v1_unnamed__12) -- (v0_unnamed__12);

  \fill[pointcolor_unnamed__12] (v1_unnamed__12) circle (0.1 pt);
  \fill[pointcolor_unnamed__12] (v0_unnamed__12) circle (0.1 pt);


  \definecolor{pointcolor_unnamed__13}{rgb}{ 1,0,0 }
  \tikzstyle{pointstyle_unnamed__13} = [fill=pointcolor_unnamed__13]

  \coordinate (v0_unnamed__13) at (-0.666667, -1);
  \coordinate (v1_unnamed__13) at (-2, -3);

  \definecolor{linecolor_unnamed__13}{rgb}{ 0 0 0 }

  \tikzstyle{linestyle_unnamed__13} = [color=linecolor_unnamed__13, thick]

  \draw[linestyle_unnamed__13] (v1_unnamed__13) -- (v0_unnamed__13);

  \fill[pointcolor_unnamed__13] (v1_unnamed__13) circle (0.1 pt);
  \fill[pointcolor_unnamed__13] (v0_unnamed__13) circle (0.1 pt);


  \definecolor{pointcolor_unnamed__14}{rgb}{ 1,0,0 }
  \tikzstyle{pointstyle_unnamed__14} = [fill=pointcolor_unnamed__14]

  \coordinate (v0_unnamed__14) at (-0.333333, -0.333333);
  \coordinate (v1_unnamed__14) at (-0.666667, -1);

  \definecolor{linecolor_unnamed__14}{rgb}{ 0 0 0 }

  \tikzstyle{linestyle_unnamed__14} = [color=linecolor_unnamed__14, thick]

  \draw[linestyle_unnamed__14] (v1_unnamed__14) -- (v0_unnamed__14);

  \fill[pointcolor_unnamed__14] (v1_unnamed__14) circle (0.1 pt);
  \fill[pointcolor_unnamed__14] (v0_unnamed__14) circle (0.1 pt);


  \definecolor{pointcolor_unnamed__15}{rgb}{ 1,0,0 }
  \tikzstyle{pointstyle_unnamed__15} = [fill=pointcolor_unnamed__15]

  \coordinate (v0_unnamed__15) at (-2, -3);
  \coordinate (v1_unnamed__15) at (-6, -7);

  \definecolor{linecolor_unnamed__15}{rgb}{ 0 0 0 }

  \tikzstyle{linestyle_unnamed__15} = [color=linecolor_unnamed__15, thick]

  \draw[linestyle_unnamed__15] (v1_unnamed__15) -- (v0_unnamed__15);

  \fill[pointcolor_unnamed__15] (v1_unnamed__15) circle (0.1 pt);
  \fill[pointcolor_unnamed__15] (v0_unnamed__15) circle (0.1 pt);


  \definecolor{pointcolor_unnamed__16}{rgb}{ 1,0,0 }
  \tikzstyle{pointstyle_unnamed__16} = [fill=pointcolor_unnamed__16]

  \coordinate (v0_unnamed__16) at (-0.333333, 0.666667);
  \coordinate (v1_unnamed__16) at (-0.333333, -0.333333);

  \definecolor{linecolor_unnamed__16}{rgb}{ 0 0 0 }

  \tikzstyle{linestyle_unnamed__16} = [color=linecolor_unnamed__16, thick]

  \draw[linestyle_unnamed__16] (v1_unnamed__16) -- (v0_unnamed__16);

  \fill[pointcolor_unnamed__16] (v1_unnamed__16) circle (0.1 pt);
  \fill[pointcolor_unnamed__16] (v0_unnamed__16) circle (0.1 pt);


  \definecolor{pointcolor_unnamed__17}{rgb}{ 1,0,0 }
  \tikzstyle{pointstyle_unnamed__17} = [fill=pointcolor_unnamed__17]

  \coordinate (v0_unnamed__17) at (6, 3);
  \coordinate (v1_unnamed__17) at (7, 3);

  \definecolor{linecolor_unnamed__17}{rgb}{ 0 0 0 }

  \tikzstyle{linestyle_unnamed__17} = [color=linecolor_unnamed__17, thick]

  \draw[linestyle_unnamed__17] (v1_unnamed__17) -- (v0_unnamed__17);

  \fill[pointcolor_unnamed__17] (v1_unnamed__17) circle (0.1 pt);
  \fill[pointcolor_unnamed__17] (v0_unnamed__17) circle (0.1 pt);


  \definecolor{pointcolor_unnamed__18}{rgb}{ 1,0,0 }
  \tikzstyle{pointstyle_unnamed__18} = [fill=pointcolor_unnamed__18]

  \coordinate (v0_unnamed__18) at (1, 0.333333);
  \coordinate (v1_unnamed__18) at (3, 1);

  \definecolor{linecolor_unnamed__18}{rgb}{ 0 0 0 }

  \tikzstyle{linestyle_unnamed__18} = [color=linecolor_unnamed__18, thick]

  \draw[linestyle_unnamed__18] (v1_unnamed__18) -- (v0_unnamed__18);

  \fill[pointcolor_unnamed__18] (v1_unnamed__18) circle (0.1 pt);
  \fill[pointcolor_unnamed__18] (v0_unnamed__18) circle (0.1 pt);


  \definecolor{pointcolor_unnamed__19}{rgb}{ 1,0,0 }
  \tikzstyle{pointstyle_unnamed__19} = [fill=pointcolor_unnamed__19]

  \coordinate (v0_unnamed__19) at (3, 1);
  \coordinate (v1_unnamed__19) at (5, 2);

  \definecolor{linecolor_unnamed__19}{rgb}{ 0 0 0 }

  \tikzstyle{linestyle_unnamed__19} = [color=linecolor_unnamed__19, thick]

  \draw[linestyle_unnamed__19] (v1_unnamed__19) -- (v0_unnamed__19);

  \fill[pointcolor_unnamed__19] (v1_unnamed__19) circle (0.1 pt);
  \fill[pointcolor_unnamed__19] (v0_unnamed__19) circle (0.1 pt);


  \definecolor{pointcolor_unnamed__20}{rgb}{ 1,0,0 }
  \tikzstyle{pointstyle_unnamed__20} = [fill=pointcolor_unnamed__20]

  \coordinate (v0_unnamed__20) at (5, 2);
  \coordinate (v1_unnamed__20) at (6, 3);

  \definecolor{linecolor_unnamed__20}{rgb}{ 0 0 0 }

  \tikzstyle{linestyle_unnamed__20} = [color=linecolor_unnamed__20, thick]

  \draw[linestyle_unnamed__20] (v1_unnamed__20) -- (v0_unnamed__20);

  \fill[pointcolor_unnamed__20] (v1_unnamed__20) circle (0.1 pt);
  \fill[pointcolor_unnamed__20] (v0_unnamed__20) circle (0.1 pt);


  \definecolor{pointcolor_unnamed__21}{rgb}{ 1,0,0 }
  \tikzstyle{pointstyle_unnamed__21} = [fill=pointcolor_unnamed__21]

  \coordinate (v0_unnamed__21) at (0, 0.333333);
  \coordinate (v1_unnamed__21) at (1, 0.333333);

  \definecolor{linecolor_unnamed__21}{rgb}{ 0 0 0 }

  \tikzstyle{linestyle_unnamed__21} = [color=linecolor_unnamed__21, thick]

  \draw[linestyle_unnamed__21] (v1_unnamed__21) -- (v0_unnamed__21);

  \fill[pointcolor_unnamed__21] (v1_unnamed__21) circle (0.1 pt);
  \fill[pointcolor_unnamed__21] (v0_unnamed__21) circle (0.1 pt);


  \definecolor{pointcolor_unnamed__22}{rgb}{ 1,0,0 }
  \tikzstyle{pointstyle_unnamed__22} = [fill=pointcolor_unnamed__22]

  \coordinate (v0_unnamed__22) at (0, 0.333333);
  \coordinate (v1_unnamed__22) at (-0.333333, 0.666667);

  \definecolor{linecolor_unnamed__22}{rgb}{ 0 0 0 }

  \tikzstyle{linestyle_unnamed__22} = [color=linecolor_unnamed__22, thick]

  \draw[linestyle_unnamed__22] (v1_unnamed__22) -- (v0_unnamed__22);

  \fill[pointcolor_unnamed__22] (v1_unnamed__22) circle (0.1 pt);
  \fill[pointcolor_unnamed__22] (v0_unnamed__22) circle (0.1 pt);


  \definecolor{pointcolor_unnamed__23}{rgb}{ 1,0,0 }
  \tikzstyle{pointstyle_unnamed__23} = [fill=pointcolor_unnamed__23]

  \coordinate (v0_unnamed__23) at (-0.333333, 0.666667);
  \coordinate (v1_unnamed__23) at (-1, 2);

  \definecolor{linecolor_unnamed__23}{rgb}{ 0 0 0 }

  \tikzstyle{linestyle_unnamed__23} = [color=linecolor_unnamed__23, thick]

  \draw[linestyle_unnamed__23] (v1_unnamed__23) -- (v0_unnamed__23);

  \fill[pointcolor_unnamed__23] (v1_unnamed__23) circle (0.1 pt);
  \fill[pointcolor_unnamed__23] (v0_unnamed__23) circle (0.1 pt);


  \definecolor{pointcolor_unnamed__24}{rgb}{ 1,0,0 }
  \tikzstyle{pointstyle_unnamed__24} = [fill=pointcolor_unnamed__24]

  \coordinate (v0_unnamed__24) at (6, 3);
  \coordinate (v1_unnamed__24) at (6, 7);

  \definecolor{linecolor_unnamed__24}{rgb}{ 0 0 0 }

  \tikzstyle{linestyle_unnamed__24} = [color=linecolor_unnamed__24, thick]

  \draw[linestyle_unnamed__24] (v1_unnamed__24) -- (v0_unnamed__24);

  \fill[pointcolor_unnamed__24] (v1_unnamed__24) circle (0.1 pt);
  \fill[pointcolor_unnamed__24] (v0_unnamed__24) circle (0.1 pt);


  \definecolor{pointcolor_unnamed__25}{rgb}{ 1,0,0 }
  \tikzstyle{pointstyle_unnamed__25} = [fill=pointcolor_unnamed__25]

  \coordinate (v0_unnamed__25) at (-1, 2);
  \coordinate (v1_unnamed__25) at (-2, 3);

  \definecolor{linecolor_unnamed__25}{rgb}{ 0 0 0 }

  \tikzstyle{linestyle_unnamed__25} = [color=linecolor_unnamed__25, thick]

  \draw[linestyle_unnamed__25] (v1_unnamed__25) -- (v0_unnamed__25);

  \fill[pointcolor_unnamed__25] (v1_unnamed__25) circle (0.1 pt);
  \fill[pointcolor_unnamed__25] (v0_unnamed__25) circle (0.1 pt);


  \definecolor{pointcolor_unnamed__26}{rgb}{ 1,0,0 }
  \tikzstyle{pointstyle_unnamed__26} = [fill=pointcolor_unnamed__26]

  \coordinate (v0_unnamed__26) at (-1, 2);
  \coordinate (v1_unnamed__26) at (-1, 7);

  \definecolor{linecolor_unnamed__26}{rgb}{ 0 0 0 }

  \tikzstyle{linestyle_unnamed__26} = [color=linecolor_unnamed__26, thick]

  \draw[linestyle_unnamed__26] (v1_unnamed__26) -- (v0_unnamed__26);

  \fill[pointcolor_unnamed__26] (v1_unnamed__26) circle (0.1 pt);
  \fill[pointcolor_unnamed__26] (v0_unnamed__26) circle (0.1 pt);


  \definecolor{pointcolor_unnamed__27}{rgb}{ 1,0,0 }
  \tikzstyle{pointstyle_unnamed__27} = [fill=pointcolor_unnamed__27]

  \coordinate (v0_unnamed__27) at (-2, 3);
  \coordinate (v1_unnamed__27) at (-3, 3);

  \definecolor{linecolor_unnamed__27}{rgb}{ 0 0 0 }

  \tikzstyle{linestyle_unnamed__27} = [color=linecolor_unnamed__27, thick]

  \draw[linestyle_unnamed__27] (v1_unnamed__27) -- (v0_unnamed__27);

  \fill[pointcolor_unnamed__27] (v1_unnamed__27) circle (0.1 pt);
  \fill[pointcolor_unnamed__27] (v0_unnamed__27) circle (0.1 pt);


  \definecolor{pointcolor_unnamed__28}{rgb}{ 1,0,0 }
  \tikzstyle{pointstyle_unnamed__28} = [fill=pointcolor_unnamed__28]

  \coordinate (v0_unnamed__28) at (-2, 3);
  \coordinate (v1_unnamed__28) at (-2, 7);

  \definecolor{linecolor_unnamed__28}{rgb}{ 0 0 0 }

  \tikzstyle{linestyle_unnamed__28} = [color=linecolor_unnamed__28, thick]

  \draw[linestyle_unnamed__28] (v1_unnamed__28) -- (v0_unnamed__28);

  \fill[pointcolor_unnamed__28] (v1_unnamed__28) circle (0.1 pt);
  \fill[pointcolor_unnamed__28] (v0_unnamed__28) circle (0.1 pt);


  \definecolor{pointcolor_unnamed__29}{rgb}{ 1,0,0 }
  \tikzstyle{pointstyle_unnamed__29} = [fill=pointcolor_unnamed__29]

  \coordinate (v0_unnamed__29) at (-3, 3);
  \coordinate (v1_unnamed__29) at (-7, -1);

  \definecolor{linecolor_unnamed__29}{rgb}{ 0 0 0 }

  \tikzstyle{linestyle_unnamed__29} = [color=linecolor_unnamed__29, thick]

  \draw[linestyle_unnamed__29] (v1_unnamed__29) -- (v0_unnamed__29);

  \fill[pointcolor_unnamed__29] (v1_unnamed__29) circle (0.1 pt);
  \fill[pointcolor_unnamed__29] (v0_unnamed__29) circle (0.1 pt);


  \definecolor{pointcolor_unnamed__30}{rgb}{ 1,0,0 }
  \tikzstyle{pointstyle_unnamed__30} = [fill=pointcolor_unnamed__30]

  \coordinate (v0_unnamed__30) at (-3, 3);
  \coordinate (v1_unnamed__30) at (-3, 7);

  \definecolor{linecolor_unnamed__30}{rgb}{ 0 0 0 }

  \tikzstyle{linestyle_unnamed__30} = [color=linecolor_unnamed__30, thick]

  \draw[linestyle_unnamed__30] (v1_unnamed__30) -- (v0_unnamed__30);

  \fill[pointcolor_unnamed__30] (v1_unnamed__30) circle (0.1 pt);
  \fill[pointcolor_unnamed__30] (v0_unnamed__30) circle (0.1 pt);


\end{tikzpicture}


%% file: smoothKleinSkeleton-A.tikz
\begin{tikzpicture}[scale=1.5]
  \tikzstyle{triangle} = [text=black,draw=black,fill=white,rectangle];

  \small

  \coordinate (v4v7v8) at (0,0);

  \coordinate (v1v4v7) at (0,1);
  \coordinate (v1v7v11) at (-1,2);

  \coordinate (v7v8v11) at (-1,-1);
  \coordinate (v8v9v11) at (-2,-3);

  \coordinate (v4v8v9) at (1,0);
  \coordinate (v1v4v9) at (3,1);

  \draw[black, ultra thick] (v1v4v7) -- (v1v7v11) -- (v7v8v11) -- (v8v9v11) -- (v4v8v9) -- (v1v4v9) -- cycle;
  \draw[black, ultra thick] (v4v7v8) -- (v1v4v7);
  \draw[black, ultra thick] (v4v7v8) -- (v7v8v11);
  \draw[black, ultra thick] (v4v7v8) -- (v4v8v9);

  \node at (v4v7v8) [triangle] {4 7 8};
  \node at (v1v4v7) [triangle] {1 4 7};
  \node at (v1v7v11) [triangle] {1 7 11};
  \node at (v7v8v11) [triangle] {7 8 11};
  \node at (v8v9v11) [triangle] {8 9 11};
  \node at (v4v8v9) [triangle] {4 8 9};
  \node at (v1v4v9) [triangle] {1 4 9};

  \node at ($ 0.5*(v4v7v8) + 0.5*(v1v4v7) $) [left] {$1/3$};
  \node at ($ 0.5*(v4v7v8) + 0.5*(v7v8v11) $) [above left] {$1/3$};
  \node at ($ 0.5*(v4v7v8) + 0.5*(v4v8v9) $) [above] {$1/3$};

  \node at ($ 0.5*(v1v4v7) + 0.5*(v1v7v11) $) [above right] {$1/3$};
  \node at ($ 0.5*(v1v7v11) + 0.5*(v7v8v11) $) [left] {$1$};

  \node at ($ 0.5*(v7v8v11) + 0.5*(v8v9v11) $) [above left] {$1/3$};
  \node at ($ 0.5*(v8v9v11) + 0.5*(v4v8v9) $) [below right] {$1$};

  \node at ($ 0.5*(v4v8v9) + 0.5*(v1v4v9) $) [below right] {$1/3$};
  \node at ($ 0.5*(v1v4v9) + 0.5*(v1v4v7) $) [above] {$1$};

\end{tikzpicture}


%% file: smoothKleinSkeleton-B.tikz
\begin{tikzpicture}[scale=2.0]
  \tikzstyle{triangle} = [text=black,draw=black,fill=white,rectangle];

  \small

  \coordinate (v4v7v8) at (0,0);
  \coordinate (v1v4v7) at (90:1);
  \coordinate (v7v8v11) at (210:1);
  \coordinate (v4v8v9) at (330:1);

  \draw[black, ultra thick] (v1v4v7) -- (v7v8v11) -- (v4v8v9) -- cycle;
  \draw[black, ultra thick] (v4v7v8) -- (v1v4v7);
  \draw[black, ultra thick] (v4v7v8) -- (v7v8v11);
  \draw[black, ultra thick] (v4v7v8) -- (v4v8v9);

  \node at (v4v7v8) [triangle] {};
  \node at (v1v4v7) [triangle] {};
  \node at (v7v8v11) [triangle] {};
  \node at (v4v8v9) [triangle] {};

  \node at ($ 0.5*(v4v7v8) + 0.5*(v1v4v7) $) [below left] {$y$}; 
  \node at ($ 0.5*(v4v7v8) + 0.5*(v7v8v11) $) [above] {$x$}; 
  \node at ($ 0.5*(v4v7v8) + 0.5*(v4v8v9) $) [above] {$z$}; 

  \node at ($ 0.5*(v1v4v7) + 0.5*(v7v8v11) $) [above left] {$u$}; 
  \node at ($ 0.5*(v7v8v11) + 0.5*(v4v8v9) $) [below] {$v$}; 
  \node at ($ 0.5*(v4v8v9) + 0.5*(v1v4v7) $) [above right] {$w$}; 

\end{tikzpicture}
